\documentclass[a4paper,10pt]{article}
\usepackage[latin1]{inputenc}
\usepackage[T1]{fontenc}
\usepackage[english]{babel}
\usepackage{amsmath}   
\usepackage{amsfonts}
\usepackage{amssymb}
\usepackage{amsthm}
\usepackage{verbatim}  
\usepackage{url}       
\usepackage{floatflt}  
\usepackage{enumerate}  
\usepackage{float}
\usepackage{color}

\newcommand{\md}{\mathrm{d}}
\newcommand{\lag}{\langle}
\newcommand{\rag}{\rangle}

\renewcommand{\AA}{\mathcal{A}}
\newcommand{\BB}{\mathcal{B}}

\newcommand{\DD}{\mathcal{D}}

\newcommand{\II}{\mathcal{I}}

\newcommand{\LL}{\mathcal{L}}

\newcommand{\OO}{\mathcal{O}}
\newcommand{\PP}{\mathcal{P}}
\newcommand{\QQ}{\mathcal{Q}}
\newcommand{\TT}{\mathcal{T}}

\newcommand{\VV}{\mathcal{V}}

\newcommand{\R}{\mathbb{R}}
\newcommand{\C}{\mathbb{C}}
\newcommand{\N}{\mathbb{N}}
\newcommand{\Z}{\mathbb{Z}}

\numberwithin{equation}{section}

\theoremstyle{definition}

\newtheorem{hypo}{Hypothesis}[section]

\theoremstyle{plain}
\newtheorem{prop}{Proposition}[section]
\newtheorem{theo}{Theorem}[section]
\newtheorem{lemme}{Lemma}[section]

\theoremstyle{remark}
\newtheorem{rmq}{Remark}[section]

\title{Simultaneous local exact controllability of 1{D} bilinear Schr\"odinger equations}
\author{
Morgan \textsc{Morancey}
\footnote{CMLS UMR 7640, Ecole Polytechnique, 91128 Palaiseau, FRANCE.  
email: Morgan.Morancey@cmla.ens-cachan.fr}
\thanks{The author was partially supported by the ``Agence Nationale de la Recherche'' (ANR),
Projet Blanc EMAQS number ANR-2011-BS01-017-01 and by CMLA UMR 8536, ENS Cachan, 61 avenue du Président Wilson, 94235 Cachan, FRANCE.} }

\begin{document}

\maketitle

\hrule
\begin{abstract}
We consider $N$ independent quantum particles, in an infinite square potential well coupled to an external laser field. These particles are modelled by a system of linear Schr\"odinger equations on a bounded interval. This is a bilinear control system in which the state is the $N$-tuple of wave functions. The control is the real amplitude of the laser field. For $N=1$, Beauchard and Laurent proved local exact controllability around the ground state in arbitrary time. We prove, under an extra generic assumption, that their result does not hold in small time if $N\geq 2$. Still, for $N=2$, we prove that local controllability holds either in arbitrary time up to a global phase or exactly up to a global delay. This is proved using Coron's return method. We also prove that for $N \geq 3$, local controllability does not hold in small time even up to a global phase. Finally, for $N=3$, we prove that local controllability holds up to a global phase and a global delay.
\end{abstract}
\hrule


\section{Introduction}

\subsection{Main results}
\label{subsect_main_result}

\noindent
We consider a quantum particle in a one dimensional infinite square potential well coupled to an external laser field. The evolution of the wave function $\psi$ is given by the following Schr\"odinger equation
\begin{equation}
\label{syst_1eq}
\left\{
\begin{aligned}
&i \partial_t \psi = - \partial^2_{xx} \psi - u(t) \mu(x) \psi, \: &(t,x)& \in (0,T)\times(0,1),
\\
&\psi(t,0) = \psi(t,1) = 0,  \quad & t & \in (0,T),
\end{aligned}
\right.
\end{equation}
where $\mu \in H^3((0,1),\R)$ is the dipolar moment and $u : t \in (0,T) \mapsto \R$ is the amplitude of the laser field. This is a bilinear control system in which the state $\psi$ lives on a sphere of $L^2((0,1),\C)$. Similar systems have been studied by various authors (see e.g. \cite{BeauchardLaurent,Mirrahimi09,RouchonModele,ChambrionSigalotti09}).

\noindent 
We are interested in simultaneous controllability of system~(\ref{syst_1eq}) and thus we consider, for $N \in \N^*$, the system
\begin{equation}
\label{syst_Neq}
\left\{
\begin{aligned}
&i \partial_t \psi^j = - \partial^2_{xx} \psi^j - u(t) \mu(x) \psi^j, \: &(t,x)& \in (0,T)\times(0,1), \, j \in \{ 1, \dots, N\},
\\
&\psi^j(t,0) = \psi^j(t,1) =0, \quad & t & \in (0,T), \, j \in \{ 1 , \dots , N \}.
\end{aligned}
\right.
\end{equation}
It is a simplified model for the evolution of $N$ identical and independent particles submitted to a single external laser field where entanglement has been neglected. This can be seen as a first step towards more sophisticated models.

\noindent 
Before going into details, let us set some notations. In this paper, $\lag \cdot , \cdot \rag$ denotes the usual scalar product on $L^2((0,1),\C)$ i.e.
\begin{equation*}
\lag f,g \rag = \int_0^1 f(x) \overline{g(x)} \md x
\end{equation*}
and $\mathcal{S}$ denotes the unit sphere of $L^2((0,1),\C)$. We consider the operator $A$ defined by
\begin{equation*}
\DD(A) := H^2\cap H^1_0((0,1),\C), \quad A \varphi := -\partial^2_{xx} \varphi.
\end{equation*}
Its eigenvalues and eigenvectors are
\begin{equation*}
\lambda_k := (k \pi)^2, \quad \varphi_k(x) := \sqrt{2} \sin(k \pi x), \quad \forall k \in \N^*.
\end{equation*}
The family $(\varphi_k)_{k \in \N^*}$ is an Hilbert basis of $L^2((0,1),\C)$. The eigenstates are defined by 
\begin{equation*}
\Phi_k(t,x) := \varphi_k(x) e^{-i \lambda_k t}, \quad (t,x) \in \R^+ \times (0,1), \,  k \in \N^*.
\end{equation*}
Any $N$-tuple of eigenstates is solution of system~(\ref{syst_Neq}) with control $u \equiv 0$. Finally, we define the spaces
\begin{equation*}
H^s_{(0)}((0,1),\C) := \DD(A^{s/2}), \; \forall s>0,
\end{equation*}
endowed with the norm
\begin{equation*}
|| \cdot ||_{H^s_{(0)}} := \left( \sum_{k=1}^{+ \infty} |k^s \lag \; \cdot \; , \varphi_k \rag |^2 \right)^{1/2}
\end{equation*}
and
\begin{equation*}
h^s(\N^*,\C) := \left\{ a = (a_k)_{k \in \N^*} \in \C^{\N^*} \, ; \, \sum_{k=1}^{+\infty} |k^s a_k|^2 < + \infty \right\}
\end{equation*}
endowed with the norm
\begin{equation*}
||a||_{h^s} := \left( \sum_{k=1}^{+ \infty} |k^s a_k|^2 \right)^{1/2}.
\end{equation*}

\noindent
Our goal is to control simultaneously the particles modelled by (\ref{syst_Neq}) with initial conditions
\begin{equation}
\label{cond_ini}
\psi^j(0,x) = \varphi_j(x), \quad x \in (0,1), \, j \in \{1,\dots , N\},
\end{equation}
locally around $\big( \Phi_1 , \dots , \Phi_N \big)$ using a single control.

\begin{rmq}
\label{rmq_invariant}
Before getting to controllability results, it has to be noticed that for any control $v \in L^2((0,T),\R)$, the associated solution of (\ref{syst_Neq}) satisfies
\begin{equation*}
\lag \psi^j(t), \psi^k(t) \rag = \lag \psi^j(0) , \psi^k(0) \rag, \quad \forall t \in [0,T].
\end{equation*}
This invariant has to be taken into account since it imposes compatibility conditions between targets and initial conditions.
\end{rmq}

\noindent
The case $N=1$ of a single equation was studied, in this setting, in \cite[Theorem 1]{BeauchardLaurent} by Beauchard and Laurent. They proved exact controllability, in $H^3_{(0)}$, in arbitrary time, locally around $\Phi_1$. Their proof relies on the linear test, the inverse mapping theorem and a regularizing effect. We prove that this result cannot be extended to the case $N=2$.

\noindent
In the spirit of \cite{BeauchardLaurent}, we assume the following hypothesis.
\begin{hypo}
\label{hypo_mu}
The dipolar moment $\mu \in H^3((0,1),\R)$ is such that there exists $c>0$ satisfying
\begin{equation*}
| \lag \mu \varphi_j , \varphi_k \rag | \geq \frac{c}{k^3}, \quad \forall k \in \N^*, \, \forall j \in \{1,\dots, N \}.
\end{equation*}
\end{hypo}

\begin{rmq}
In the same way as in \cite[Proposition 16]{BeauchardLaurent}, one may prove that Hypothesis~\ref{hypo_mu} holds generically in $H^3((0,1),\R)$.
\end{rmq} 

\noindent
Using \cite[Theorem 1]{BeauchardLaurent}, Hypothesis~\ref{hypo_mu} implies that the $j^{th}$ equation of system~(\ref{syst_Neq}) is locally controllable in $H^3_{(0)}$ around $\Phi_j$.

\begin{hypo}
\label{hypo_mu2}
The dipolar moment $\mu \in H^3((0,1),\R)$ is such that
\begin{equation*}
\AA := \lag \mu \varphi_1 , \varphi_1 \rag \lag (\mu')^2 \varphi_2 , \varphi_2 \rag - \lag \mu \varphi_2 , \varphi_2 \rag \lag (\mu')^2 \varphi_1 , \varphi_1 \rag  \neq 0.
\end{equation*}
\end{hypo}

\begin{rmq}
For example, $\mu(x):=x^3$ satisfies both Hypothesis~\ref{hypo_mu} and \ref{hypo_mu2}. Unfortunately, the case $\mu(x):=x$ studied in \cite{RouchonModele} does not satisfy these hypotheses. But, as in \cite[Proposition 16]{BeauchardLaurent}, one may prove that Hypotheses \ref{hypo_mu} and \ref{hypo_mu2} hold simultaneously generically in $H^3((0,1),\R)$.
\end{rmq}

\begin{rmq}
Hypothesis~\ref{hypo_mu2} implies that there exists $j \in \{1,2\}$ such that $\lag \mu \varphi_j , \varphi_j \rag \neq 0$. Without loss of generality, when Hypothesis~\ref{hypo_mu2} is assumed to hold, one should consider that $\lag \mu \varphi_1, \varphi_1 \rag \neq 0$.
\end{rmq}

\begin{theo}
\label{th_NC_2eq_1}
Let $N=2$ and $\mu \in H^3((0,1),\R)$ be such that Hypothesis~\ref{hypo_mu2} hold. Let $\alpha \in \{-1,1\}$ be defined by $\alpha := \text{sign} (\AA \lag \mu \varphi_1, \varphi_1\rag)$.  There exists $T_*>0$ and $\varepsilon>0$ such that for any $T<T_*$, for every $u \in L^2((0,T),\R)$ with $ ||u||_{L^2(0,T)} < \varepsilon$,
the solution of system~(\ref{syst_Neq})-(\ref{cond_ini}) satisfies
\begin{equation*}
\big( \psi^1(T), \psi^2(T) \big) \neq \left( \Phi_1(T) ,
\left( \sqrt{1-\delta^2} + i \alpha \delta \right) \Phi_2(T) \right), \, \forall \delta >0 .
\end{equation*}
\end{theo}

\noindent
Thus, under Hypothesis~\ref{hypo_mu2}, simultaneous controllability does not hold for $(\psi^1,\psi^2)$ around $(\Phi_1,\Phi_2)$ in small time with small controls. The smallness assumption on the control is in $L^2$ norm. This prevents from extending \cite[Theorem 1]{BeauchardLaurent} to the case $N \geq 2$. Notice that the proposed target that cannot be reached satisfies the compatibility conditions of Remark~\ref{rmq_invariant}.

\noindent
However, when modelling a quantum particle, the global phase is physically meaningless. Thus for any $\theta \in \R$ and $\psi^1,\psi^2 \in L^2((0,1),\C)$, the states $e^{i \theta} (\psi^1, \psi^2)$ and $(\psi^1,\psi^2)$ are physically equivalent. Working up to a global phase, we prove the following theorem.
\begin{theo}
\label{th_controle_2eq_1}
Let $N=2$. Let $T>0$. Let $\mu \in H^3((0,1),\C)$ satisfy Hypothesis~\ref{hypo_mu} and $\lag \mu \varphi_1, \varphi_1 \rag \neq \lag \mu \varphi_2, \varphi_2 \rag$. There exists $\theta \in \R$, $\varepsilon_0>0$ and a $C^1$ map 
\begin{equation*}
\Gamma : \OO_{\varepsilon_0} \to L^2((0,T),\R)
\end{equation*}
where 
\begin{align*}
\OO_{\varepsilon_0} := \Big\{ \big( \psi^1_f , \psi^2_f \big) \in H^3_{(0)}((0,1),\C)^2 \, ; \, \lag \psi^j_f , \psi^k_f \rag = \delta_{j=k} \text{ and }
\\
 \sum_{j=1}^2 ||\psi^j_f - e^{i \theta} \Phi_j(T)||_{H^3_{(0)}} < \varepsilon_0 \Big\},
\end{align*}
such that for any $\big( \psi^1_f,\psi^2_f \big) \in \OO_{\varepsilon_0}$, the solution of system~(\ref{syst_Neq}) with initial condition (\ref{cond_ini}) and control $u= \Gamma\big( \psi^1_f,\psi^2_f \big)$ satisfies
\begin{equation*}
(\psi^1(T),\psi^2(T)) = (\psi^1_f, \psi^2_f).
\end{equation*}
\end{theo}

\begin{rmq} \label{rmq_theo}
Notice that, using Remark~\ref{rmq_invariant}, the condition $\lag \psi^j_f, \psi^k_f \rag = \delta_{j=k}$ is not restrictive. Indeed, as $\psi^j(0) = \varphi_j$, we can only reach targets satisfying such an orthonormality condition.
\end{rmq}

\begin{rmq}
\label{rmq_gestion_cond_ini}
The same theorem holds with initial conditions $(\psi^1_0,\psi^2_0)$ close enough to $(\varphi_1, \varphi_2)$ in $H^3_{(0)}$ satisfying the constraints $\lag \psi^1_0, \psi^2_0\rag = \lag \psi^1_f, \psi^2_f \rag$ (see Remark~\ref{rmq_th_cond_ini} in Section \ref{subsect_controle_NL}).
\end{rmq}

\noindent
Working in time large enough we can drop the global phase and prove the following theorem.
\begin{theo}
\label{th_controle_2eq_2}
Let $N=2$. Let $\mu \in H^3((0,1),\C)$ satisfy Hypothesis~\ref{hypo_mu} and $4 \lag \mu \varphi_1 ,\varphi_1 \rag - \lag \mu \varphi_2 , \varphi_2 \rag \neq 0$. There exists $T^*>0$ such that, for any $T \geq 0$, there exists $\varepsilon_0 > 0$ and a $C^1$ map 
\begin{equation*}
\Gamma : \OO_{\varepsilon_0,T} \to L^2((0,T^*+T),\R)
\end{equation*}
where 
\begin{align*}
\OO_{\varepsilon_0,T} := \Big\{ \big( \psi^1_f , \psi^2_f \big) \in H^3_{(0)}((0,1),\C)^2 \, ; \, \lag \psi^j_f , \psi^k_f \rag = \delta_{j=k} \text{ and } 
\\
\sum_{j=1}^2 ||\psi^j_f - \Phi_j(T) ||_{H^3_{(0)}} < \varepsilon_0 \Big\},
\end{align*}
such that for any $\big( \psi^1_f,\psi^2_f \big) \in \OO_{\varepsilon_0,T}$, the solution of system~(\ref{syst_Neq}) with initial condition (\ref{cond_ini}) and control $u= \Gamma\big( \psi^1_f,\psi^2_f \big)$ satisfies
\begin{equation*}
\big( \psi^1(T^*+T) , \psi^2(T^*+T) \big) = \big( \psi^1_f , \psi^2_f \big).
\end{equation*}
\end{theo}

\begin{rmq}
Remark~\ref{rmq_gestion_cond_ini} is still valid in this case.
\end{rmq}

\noindent
We now turn to the case $N=3$. We prove that under an extra generic assumption, Theorem~\ref{th_controle_2eq_1} cannot be extended to three particles.
Assume the following hypothesis.
\begin{hypo}
\label{hypo_mu3}
The dipolar moment $\mu \in H^3((0,1),\R)$ is such that
\begin{align*}
\BB :&= \big( \lag \mu \varphi_3, \varphi_3 \rag - \lag \mu \varphi_2, \varphi_2 \rag \big)
\lag (\mu')^2 \varphi_1, \varphi_1 \rag
\\
&+ \big( \lag \mu \varphi_1, \varphi_1 \rag - \lag \mu \varphi_3, \varphi_3 \rag \big)
\lag (\mu')^2 \varphi_2, \varphi_2 \rag
\\
&+ \big( \lag \mu \varphi_2, \varphi_2 \rag - \lag \mu \varphi_1, \varphi_1 \rag \big)
\lag (\mu')^2 \varphi_3, \varphi_3 \rag
\neq 0.
\end{align*}
\end{hypo}

\begin{rmq}
Hypothesis~\ref{hypo_mu3} implies that there exist $j , k \in \{1,2,3\}$ such that $\lag \mu \varphi_j , \varphi_j \rag \neq \lag \mu \varphi_k, \varphi_k \rag$. Without loss of generality, when Hypothesis~\ref{hypo_mu3} is assumed to hold, one should consider that $\lag \mu \varphi_1, \varphi_1 \rag \neq \lag \mu \varphi_2, \varphi_2 \rag$.
\end{rmq}

\begin{rmq}
Again, Hypotheses \ref{hypo_mu} and \ref{hypo_mu3} hold simultaneously generically in $H^3((0,1),\R)$.
\end{rmq}

\noindent
We prove the following theorem.
\begin{theo}
\label{th_NC_3eq}
Let $N=3$ and $\mu \in H^3((0,1),\R)$ be such that Hypothesis~\ref{hypo_mu3} hold. Let $\beta \in \{-1 ,1\}$ be defined by $\beta = \text{sign} \big(\BB (\lag \mu \varphi_2,\varphi_2 \rag - \lag \mu \varphi_1, \varphi_1\rag)\big)$. There exists $T_*>0$ and $\varepsilon>0$ such that, for any $T<T_*$, for every $u \in L^2((0,T),\R)$ with $||u||_{L^2(0,T)} < \varepsilon$,
the solution of system~(\ref{syst_Neq})-(\ref{cond_ini}) satisfies for every $\delta >0$ and $\nu \in \R$,\begin{equation*}
\big( \psi^1(T), \psi^2(T) ,\psi^3(T) \big) \neq  e^{i \nu} \left( \Phi_1(T) , \Phi_2(T), 
\left( \sqrt{1-\delta^2} + i \beta \delta \right) \Phi_3(T) \right).
\end{equation*}
\end{theo}
Thus, in small time, local exact controllability with small controls does not hold for $N\geq 3$, even up to a global phase. The next statement ensures that it holds up to a global phase and a global delay.
\begin{theo}
\label{theo_controle_3eq}
Let $N=3$. Let $\mu \in H^3((0,1),\C)$ satisfy Hypothesis~\ref{hypo_mu} and $5 \lag \mu \varphi_1, \varphi_1\rag - 8 \lag \mu \varphi_2, \varphi_2 \rag + 3 \lag \mu \varphi_3 , \varphi_3 \rag \neq 0$. There exists $\theta \in \R$, $T^*>0$ such that, for any $T\geq 0$, there exists $\varepsilon_0 > 0$ and a $C^1$ map 
\begin{equation*}
\Gamma : \OO_{\varepsilon_0,T} \to L^2((0,T^*+T),\R)
\end{equation*}
where 
\begin{align*}
\OO_{\varepsilon_0,T} := \Big\{ \big( \psi^1_f , \psi^2_f , \psi^3_f \big) \in H^3_{(0)}((0,1),\C)^3 \, ; \, \lag \psi^j_f , \psi^k_f \rag = \delta_{j=k} \text{ and }
\\
 \sum_{j=1}^3 ||\psi^j_f - e^{i \theta} \Phi_j(T)||_{H^3_{(0)}} < \varepsilon_0 \Big\},
\end{align*}
such that for any $\big( \psi^1_f,\psi^2_f,\psi^3_f \big) \in \OO_{\varepsilon_0,T}$, the solution of system~(\ref{syst_Neq}) with initial condition (\ref{cond_ini}) and control $u= \Gamma\big( \psi^1_f,\psi^2_f,\psi^3_f \big)$ satisfies 
\begin{equation*}
\big( \psi^1(T^*+T), \psi^2(T^*+T) ,\psi^3(T^*+T) \big)= \big( \psi^1_f , \psi^2_f ,\psi^3_f \big).
\end{equation*}
\end{theo}

\begin{rmq}
Remark~\ref{rmq_gestion_cond_ini} is still valid in this case.
\end{rmq}

\subsection{Heuristic}

\noindent
Contrarily to the case $N=1$, the linearized system around a $N$-tuple of eigenstates is not controllable when $N\geq 2$. Let us consider, for $N=2$, the linearization of system~(\ref{syst_Neq}) around $(\Phi_1,\Phi_2)$
\begin{equation} \label{linearise_u=0}
\left\{
\begin{aligned}
&i \partial_t \Psi^j = - \partial^2_{xx} \Psi^j - v(t) \mu(x) \Phi_j, \quad &(t,x)& \in (0,T) \times (0,1), \,
 j \in \{1,2\},
\\
&\Psi^j(t,0) = \Psi^j(t,1) = 0, \quad &t& \in (0,T), 
\\
&\Psi^j(0,x) =  0, \quad & x &\in (0,1).
\end{aligned}
\right.
\end{equation}
For $j=1,2$, straightforward computations lead to
\begin{equation} \label{ordre1_heuristique}
\Psi^j(T) = i \sum_{k=1}^{+ \infty} \lag \mu \varphi_j , \varphi_k \rag \int_0^T v(t) e^{i(\lambda_k - \lambda_j)t} \md t \Phi_k(T).
\end{equation}
Thus, thanks to Hypothesis~\ref{hypo_mu}, we could, by solving a suitable moment problem, control any direction $\lag \Psi^j(T) ,\Phi_k(T) \rag$, for $k \geq 2$ (with a slight abuse of notation for the direction $\Phi_k$ of the $j^{th}$ equation). Straightforward computations using (\ref{ordre1_heuristique}) lead to
\begin{equation*}
\lag \Psi^1(T) , \Phi_2(T) \rag + \overline{\lag \Psi^2(T) , \Phi_1(T) \rag} = 0.
\end{equation*}
This comes from the linearization of the invariant (see Remark~\ref{rmq_invariant})
\begin{equation*}
\lag \psi^1(t) , \psi^2(t) \rag = \lag \psi^1_0 , \psi^2_0 \rag, \quad \forall t \in (0,T),
\end{equation*}
and can be overcome (see Subsection \ref{subsect_controle_NL}). However, (\ref{ordre1_heuristique}) also implies that
\begin{equation*}
\lag \mu \varphi_2, \varphi_2 \rag \lag \Psi^1(T) , \Phi_1(T) \rag  = \lag \mu \varphi_1, \varphi_1 \rag \lag \Psi^2(T), \Phi_2(T) \rag,
\end{equation*}
for any $v \in L^2((0,T),\R)$. This is a strong obstacle to controllability and leads to Theorem~\ref{th_NC_2eq_1} (see Section \ref{sect_non_controlabilite}).

\noindent
In this situation, where a direction is lost at the first order, one can try to recover it at the second order. This strategy was used for example by Cerpa and Crépeau in \cite{CerpaCrepeau} on a Korteweg De Vries equation and adapted on the considered bilinear Schr\"odinger equation (\ref{syst_1eq}) by Beauchard and the author in \cite{BM_Tmin}. Let, for $j \in \{1,2\}$,
\begin{equation*}
\left\{
\begin{aligned}
&i \partial_t \xi^j = - \partial^2_{xx} \xi^j - v(t) \mu(x) \Psi^j - w(t) \mu(x) \Phi_j, \, &(t,x)& \in (0,T) \times (0,1), 
\\
&\xi^j(t,0) = \xi^j(t,1) = 0, \, &t& \in (0,T), 
\\
&\xi^j(0,x) =  0, \, & x &\in (0,1).
\end{aligned}
\right.
\end{equation*}
The main idea of this strategy is to exploit a rotation phenomenon when the control is turned off. However, as proved in \cite[Lemma 4]{BM_Tmin}, there is no rotation phenomenon on the diagonal directions $\lag \xi^j(T) , \Phi_j(T) \rag$ and this power series expansion strategy cannot be applied to this situation.

\noindent
Thus, the local exact controllability results in this article are proved using Coron's return method. This strategy, detailed in \cite[Chapter 6]{CoronBook}, relies on finding a reference trajectory of the non linear control system with suitable origin and final positions such that the linearized system around this reference trajectory is controllable. Then, the inverse mapping theorem allows to prove local exact controllability.

\noindent
As the Schr\"odinger equation is not time reversible, the design of the reference trajectory $(\psi^1_{ref} , \dots , \psi^N_{ref},u_{ref})$ is not straightforward. The reference control $u_{ref}$ is designed in two steps. The first step is to impose restrictive conditions on $u_{ref}$ on an arbitrary time interval $(0,\varepsilon)$ in order to ensure the controllability of the linearized system. Then, $u_{ref}$ is designed on $(\varepsilon,T^*)$ such that the reference trajectory at the final time coincides with the target. For example, to prove Theorem~\ref{theo_controle_3eq}, the reference trajectory is designed such that
\begin{equation}
\label{traj_ref_temps_final}
\big( \psi^1_{ref}(T^*), \psi^2_{ref}(T^*), \psi^3_{ref}(T^*) \big) = e^{i \theta} \big( \varphi_1 , \varphi_2 , \varphi_3 \big).
\end{equation}

\subsection{Structure of the article}

\noindent
This article is organized as follows. We recall, in Section \ref{sect_bien_pose}, well posedness results.

\noindent
To emphasize the ideas developed in this article, we start by proving Theorem~\ref{theo_controle_3eq}.
Section \ref{sect_construction_uref} is devoted to the construction of the reference trajectory.
In Subsection \ref{subsect_controle_L}, we prove the controllability of the linearized system around the reference trajectory. 
In Subsection \ref{subsect_controle_NL}, we conclude the return method thanks to an inverse mapping argument.

\noindent
In Section \ref{sect_2eq}, we adapt the construction of the reference trajectory for two equations leading to Theorems \ref{th_controle_2eq_1} and \ref{th_controle_2eq_2}.

\noindent
Finally, Section \ref{sect_non_controlabilite} is devoted to non controllability results and the proofs of Theorems \ref{th_NC_2eq_1} and \ref{th_NC_3eq}.

\subsection{A review of previous results}

Let us recall some previous results about the controllability of Schr\"odinger equations. In \cite{BallMarsdenSlemrod82}, Ball, Marsden and Slemrod proved a negative result for infinite dimensional bilinear control systems. The adaptation of this result to Schr\"odinger equations, by Turinici \cite{Turinici00}, proves that the reachable set with $L^2$ controls has an empty interior in $\mathcal{S} \cap H^2_{(0)}((0,1),\C)$. Although this is a negative result it does not prevent controllability in more regular spaces.

Actually, in \cite{Beauchard05}, Beauchard proved local exact controllability in $H^7$ using Nash-Moser theorem for a one dimensional model. The proof of this result was simplified, by Beauchard and Laurent in \cite{BeauchardLaurent}, by exhibiting a regularizing effect allowing to apply the classical inverse mapping theorem. In \cite{BeauchardCoron06}, Beauchard and Coron also proved exact controllability between eigenstates for a particle in a moving potential well.

Using stabilization techniques and Lyapunov functions, Nersesyan proved in \cite{Nersesyan10} that Beauchard and Laurent's result holds globally in $H^{3+\varepsilon}$. Other stabilization results on similar models were obtained in \cite{BeauchardMirrahimi09,Mirrahimi09,Nersesyan,BeauchardNersesyan,Morancey_polarisabilite} by Mirrahimi, Beauchard, Nersesyan and the author.

Unlike exact controllability, approximate controllability results have been obtained for Schr\"odinger equations on multidimensional domains. In \cite{CMSB09}, Chambrion, Mason, Sigalotti and Boscain proved approximate controllability in $L^2$, thanks to geometric technics on the Galerkin approximation both for the wave function and density matrices. These results were extended to stronger norms in \cite{BoussaidCaponigroChambrion} by Boussaid, Caponigro and Chambrion. Approximate controllability in more regular spaces (containing $H^3$) were obtained by Nersesyan and Nersisyan~\cite{NersesyanNersisyan12} using exact controllability in infinite time. Approximate controllability has also been obtained by Ervedoza and Puel in \cite{ErvedozaPuel09} on a model of trapped ions.


Simultaneous exact controllability of quantum particles has been obtained on a finite dimensional model in \cite{TuriniciRabitz04} by Turinici and Rabitz. Their model uses specific orientation of the molecules and their proof relies on iterated Lie brackets. In addition to the results of \cite{CMSB09}, simultaneous approximate controllability was also studied in \cite{ChambrionSigalotti09} by Chambrion and Sigalotti. They used controllability of the Galerkin approximations for a model of different particles with the same control operator and a model of identical particles with different control operators. These simultaneous approximate controllability results are valid regardless of the number of particles considered.

Finally, let us give some details about the return method. This idea of designing a reference trajectory such that the linearized system is controllable was developed by Coron in \cite{Coron92} for a stabilization problem. It was then successfully used to prove exact controllability for various systems : Euler equations in \cite{Coron93,Glass00,Glass07} by Coron and Glass, Navier-Stokes equations in \cite{Coron96, FursikovImanuvilov99, Chapouly09, CoronGuerrero09} by Coron, Fursikov, Imanuvilov, Chapouly and Guerrero, Bürgers equations in \cite{Horsin98,  GlassGuerrero07, ChapoulyBurgers09} by Horsin, Glass, Guerrero and Chapouly and many other models such as \cite{Coron02, Glass03, Glass08, CoronGuerreroRosier10}. This method was also used for a bilinear Schr\"odinger equation in \cite{Beauchard05} by Beauchard.

The question of simultaneous exact controllability for linear PDE is already present in the book \cite{LionsBook1} by Lions. He considered the case of two wave equations with different boundary controls. This was later extended to other systems by Avdonin, Tucsnak, Moran and Kapitonov in \cite{AvdoninTucsnak01, AvdoninMoran01, Kapitonov99}.

To conclude, the question of impossibility of certain motions in small time, at stake in this article, for bilinear Schr\"odinger equations was studied in \cite{Coron06,BM_Tmin} by Coron, Beauchard and the author.




\section{Well posedness}
\label{sect_bien_pose}

First, we recall the well posedness of the considered Schr\"odinger equation with a source term which proof is in \cite[Proposition 2]{BeauchardLaurent}. Consider
\begin{equation}
\label{schrodinger_f}
\left\lbrace
\begin{aligned}
&i \partial_t \psi(t,x) = - \partial_{xx}^2 \psi(t,x) - u(t) \mu(x) \psi(t,x) - f(t,x), \,  &(t,x)& \in (0,T)\times (0,1),
\\
&\psi(t,0) = \psi(t,1) = 0 , \quad &t& \in (0,T),
\\
&\psi(0,x) = \psi_0(x), \quad &x& \in (0,1).
\end{aligned}
\right.
\end{equation}

\begin{prop}
\label{bien_pose}
Let $\mu \in H^3((0,1),\R)$, $T>0$, $\psi_0 \in H^3_{(0)}(0,1)$, $u \in L^2((0,T),\R)$ and $f \in L^2((0,T), H^3 \cap H^1_0)$.
There exists a unique weak solution of (\ref{schrodinger_f}), i.e.
a function $\psi \in C^0([0,T],H^3_{(0)})$ such that the following equality 
holds in $H^{3}_{(0)}((0,1),\C)$ for every $t \in [0,T]$,
\begin{equation*}
\psi(t)=e^{-iAt} \psi_{0} +i \int_{0}^{t} e^{-iA(t-\tau)}
[ u(\tau) \mu \psi(\tau) + f(\tau) ] d\tau.
\end{equation*}
Moreover, for every $R>0$, there exists $C=C(T,\mu,R)>0$ such that,
if $\|u\|_{L^2(0,T)} < R$, then this weak solution satisfies
\begin{equation*}
\| \psi \|_{C^{0}([0,T],H^{3}_{(0)})} \leqslant 
C \Big( \| \psi_{0} \|_{H^{3}_{(0)}} + \|f\|_{L^{2}((0,T),H^{3} \cap H^1_0)} \Big).
\end{equation*}
If $f \equiv 0$, then 
\begin{equation*}
\|\psi(t)\|_{L^2(0,1)} = \|\psi_0\|_{L^2(0,1)}, \quad \forall t \in [0,T].
\end{equation*}
\end{prop}

\section{Construction of the reference trajectory for three equations}
\label{sect_construction_uref}

\noindent
The goal of this section is the design of the following family of reference trajectories to prove Theorem~\ref{theo_controle_3eq}.
\begin{theo}
\label{th_construction_uref}
Let $N=3$. Let $\mu \in H^3((0,1),\C)$ satisfy Hypothesis~\ref{hypo_mu} and $5 \lag \mu \varphi_1, \varphi_1\rag - 8 \lag \mu \varphi_2, \varphi_2 \rag + 3 \lag \mu \varphi_3 , \varphi_3 \rag \neq 0$. Let $T_1>0$ be arbitrary, $\varepsilon \in (0,T_1)$ and $\varepsilon_1 \in (\frac{\varepsilon}{2} , \varepsilon)$. There exist $\overline{\eta}>0$, $C>0$ such that for every $\eta \in (0,\overline{\eta})$, there exist $T^{\eta}>T_1$, $\theta^\eta \in \R$ and $u_{ref}^{\eta} \in L^2((0,T^{\eta}),\R)$ with 
\begin{equation} \label{borne_uref}
||u_{ref}^{\eta}||_{L^2(0,T^{\eta})}\leq C \eta
\end{equation}
such that the associated solution $\big( \psi^{1,\eta}_{ref},\psi^{2,\eta}_{ref}, \psi^{3,\eta}_{ref} \big)$ of (\ref{syst_Neq})-(\ref{cond_ini}) satisfies
\begin{equation}
\label{uref_condition1}
\begin{aligned}
\lag \mu \psi^{1,\eta}_{ref} (\varepsilon_1) , \psi^{1,\eta}_{ref} (\varepsilon_1) \rag &= \lag \mu \varphi_1 , \varphi_1 \rag + \eta,
\\
\lag \mu \psi^{2,\eta}_{ref} (\varepsilon_1) , \psi^{2,\eta}_{ref} (\varepsilon_1) \rag &= \lag \mu \varphi_2 , \varphi_2 \rag,
\\
\lag \mu \psi^{3,\eta}_{ref} (\varepsilon_1) , \psi^{3,\eta}_{ref} (\varepsilon_1) \rag &= \lag \mu \varphi_3 , \varphi_3 \rag,
\end{aligned}
\end{equation}

\begin{equation}
\label{uref_condition2}
\begin{aligned}
\lag \mu \psi^{1,\eta}_{ref} (\varepsilon) , \psi^{1,\eta}_{ref} (\varepsilon) \rag &= \lag \mu \varphi_1 , \varphi_1 \rag,
\\
\lag \mu \psi^{2,\eta}_{ref} (\varepsilon) , \psi^{2,\eta}_{ref} (\varepsilon) \rag &= \lag \mu \varphi_2 , \varphi_2 \rag + \eta,
\\
\lag \mu \psi^{3,\eta}_{ref} (\varepsilon) , \psi^{3,\eta}_{ref} (\varepsilon) \rag &= \lag \mu \varphi_3 , \varphi_3 \rag,
\end{aligned}
\end{equation}
and
\begin{equation}
\label{3eq_temps_final}
\big(\psi^{1,\eta}_{ref}(T^{\eta}) , \psi^{2,\eta}_{ref}(T^{\eta}) , \psi^{3,\eta}_{ref}(T^{\eta}) \big) = e^{i \theta^\eta} \big( \varphi_1, \varphi_2, \varphi_3 \big).
\end{equation}
\end{theo}

\begin{rmq}
For any $T\geq 0$, $u_{ref}^{\eta}$ is extended by zero on $(T^{\eta},T^{\eta}+T)$. Thus, there exists $C>0$ such that, $||u_{ref}^{\eta}||_{L^2(0,T^{\eta}+T)} \leq C \eta$, (\ref{uref_condition1}), (\ref{uref_condition2}) are satisfied and
\begin{equation*}
\big( \psi^{1,\eta}_{ref}(T^{\eta}+T) , \psi^{2,\eta}_{ref}(T^{\eta}+T) , \psi^{3,\eta}_{ref}(T^{\eta}+T) \big) = e^{i \theta^\eta} \big( \Phi_1(T) , \Phi_2(T) , \Phi_3(T) \big).
\end{equation*}
\end{rmq}

\begin{rmq}
The choice of a parameter $\eta$ sufficiently small together with conditions (\ref{uref_condition1}) and (\ref{uref_condition2}) will be used in Section \ref{subsect_controle_L} to prove the controllability of the linearized system around the reference trajectory.
The control $u_{ref}^{\eta}$ will be designed on $(0,T_1)$ and extended by zero on $(T_1,T^{\eta})$.
\end{rmq}

\noindent
The proof of Theorem~\ref{th_construction_uref} is divided in two steps : the construction of $u_{ref}^{\eta}$ on $(0,\varepsilon)$ to prove (\ref{uref_condition1}) and (\ref{uref_condition2}) and then, the construction on $(\varepsilon,T_1)$ to prove~(\ref{3eq_temps_final}). This is what is detailed in the next subsections.

\subsection{Construction on $(0,\varepsilon)$}

Let $u_{ref}^{\eta} \equiv 0$ on $[0, \frac{\varepsilon}{2})$. We prove the following proposition.
\begin{prop}
\label{conditions_minimalite}
Let $\mu \in H^3((0,1),\C)$ satisfy Hypothesis~\ref{hypo_mu}.
There exists $\eta^* >0$ and a $C^1$ map
\begin{equation*}
\hat{\Gamma} : (0,\eta^*) \rightarrow L^2 \left( \left(\frac{\varepsilon}{2},\varepsilon \right),\R \right),
\end{equation*}
such that $\hat{\Gamma}(0)= 0$ and for any $\eta \in (0,\eta^*)$, the solution $\big(\psi^{1,\eta}_{ref}, \psi^{2,\eta}_{ref} ,\psi^{3,\eta}_{ref}\big)$ of system~(\ref{syst_Neq}) with control $u_{ref}^{\eta}:= \hat{\Gamma}(\eta)$ and initial conditions $\psi^{j,\eta}_{ref}(\frac{\varepsilon}{2}) = \Phi_j(\frac{\varepsilon}{2})$, for $j=1,2,3$, satisfies (\ref{uref_condition1}) and (\ref{uref_condition2}).
\end{prop}

\noindent
\textbf{Proof of Proposition~\ref{conditions_minimalite} :}
Using Proposition~\ref{bien_pose}, it comes that the map
\begin{equation*}
\begin{array}{cccc}
\tilde{\Theta}: & L^2((\frac{\varepsilon}{2},\varepsilon),\R) 
& \rightarrow & \R^3 \times \R^3
\\
& u & \mapsto & \left(  \tilde{\Theta}_1(u) , \tilde{\Theta}_2(u) \right)
\end{array}
\end{equation*}
where
\begin{equation*}
\tilde{\Theta}_1(u) := 
\left(
\lag \mu \psi^j (\varepsilon_1) , \psi^j (\varepsilon_1) \rag - \lag \mu \varphi_j , \varphi_j \rag
\right)_{j=1,2,3} ,
\end{equation*}
and
\begin{equation*}
\tilde{\Theta}_2(u) := 
\left(
\lag \mu \psi^j (\varepsilon) , \psi^j (\varepsilon) \rag - \lag \mu \varphi_j , \varphi_j \rag
\right)_{j=1,2,3} ,
\end{equation*}
is well defined, $C^1$, satisfies $\tilde{\Theta}(0)= 0$ and 
\begin{equation}
\label{minimalite_differentielle}
\md \tilde{\Theta} (0) .v = \left( 
 \left( 2 \Re ( \lag \mu \Psi^j(\varepsilon_1) , \Phi_j(\varepsilon_1) \rag ) \right)_{1\leq j \leq 3} ,
 \left( 2 \Re ( \lag \mu \Psi^j(\varepsilon) , \Phi_j(\varepsilon) \rag ) \right)_{1, \leq j \leq 3} 
\right),
\end{equation}
where $\big( \Psi^1, \Psi^2, \Psi^3 \big)$ is the solution of (\ref{linearise_u=0})  on the time interval $\left(\frac{\varepsilon}{2},\varepsilon \right)$ with control $v$ and initial conditions $\Psi^j \left( \frac{\varepsilon}{2},\cdot \right) = 0$.
Let us prove that $\md \tilde{\Theta}(0)$ is surjective; then the inverse mapping theorem will give the conclusion.
\\

\noindent
Let $\gamma = (\gamma_j)_{1\leq j \leq 6} \in  \R^6$ and $K \geq 4$. By Proposition~\ref{prop_pb_moment} (see the appendix), there exist $v_1 \in L^2((\frac{\varepsilon}{2},\varepsilon_1),\R)$ and $v_2 \in L^2((\varepsilon_1, \varepsilon),\R) $ such that
\begin{align*}
\int_{\frac{\varepsilon}{2}}^{\varepsilon_1} v_1(t) e^{i(\lambda_k - \lambda_j)t} \md t 
&= 0, \quad \forall k \in \N^* \backslash \{ K \}, \, \forall 1 \leq j \leq 3,
\\
\int_{\frac{\varepsilon}{2}}^{\varepsilon_1} v_1(t) e^{i(\lambda_K - \lambda_j)t} \md t 
&= \frac{e^{i(\lambda_K-\lambda_j) \varepsilon_1} \gamma_j}{ 2 i \lag \mu \varphi_j , \varphi_K \rag^2}, \quad
\forall 1 \leq j \leq 3,
\\
\int_{\varepsilon_1}^{\varepsilon} v_2(t) e^{i(\lambda_k - \lambda_j)t} \md t 
&= 0, \quad \forall k \in \N^* \backslash \{ K \}, \, \forall 1 \leq j \leq 3,
\\
\int_{\varepsilon_1}^{\varepsilon} v_2(t) e^{i(\lambda_K - \lambda_j)t} \md t 
&= \frac{e^{i(\lambda_K - \lambda_j) \varepsilon} \gamma_{3+j}}{2 i \lag \mu \varphi_j , \varphi_K \rag^2} - \frac{e^{i(\lambda_K - \lambda_j) \varepsilon_1} \gamma_j}{2 i \lag \mu \varphi_j , \varphi_K \rag^2} , \quad \forall 1 \leq j \leq 3.
\end{align*}
Notice that the moments associated to redundant frequencies in the previous moment problem are all set to the same value and, as $K \geq 4$, the frequencies $\lambda_K-\lambda_j$ for $1 \leq j \leq 3$ are distinct.
Let $v \in L^2 \left(\frac{\varepsilon}{2}, \varepsilon \right)$ be defined by $v_1$ on $\left( \frac{\varepsilon}{2}, \varepsilon_1 \right)$ and by $v_2$ on $(\varepsilon_1, \varepsilon)$. Straightforward computations lead to $\md \tilde{\Theta}(0).v = \gamma$.

\hfill $\blacksquare$

\noindent

\subsection{Construction on $(\varepsilon,T_1)$}

\noindent
For any $j \in \N^*$, let $\PP_j$ be the orthogonal projection of $L^2((0,1),\C)$ onto $\text{Span}_{\C}(\varphi_k ,k\geq j+1)$ i.e.
\begin{equation*}
\PP_j(\psi) := \sum_{k=j+1}^{+ \infty} \lag \psi , \varphi_k \rag \varphi_k.
\end{equation*}
The goal of this subsection is the proof of the following proposition.
\begin{prop}
\label{traj_ref_inversion_locale}
Let $0< T_0 < T_f$. Let $\mu \in H^3((0,1),\C)$ satisfy Hypothesis~\ref{hypo_mu} and $5 \lag \mu \varphi_1, \varphi_1\rag - 8 \lag \mu \varphi_2, \varphi_2 \rag + 3 \lag \mu \varphi_3 , \varphi_3 \rag \neq 0$. There exist $\delta > 0$ and a $C^1$-map 
\begin{equation*}
\tilde{\Gamma}_{T_0,T_f} : \tilde{\OO}_{\delta,T_0} \to L^2((T_0,T_f) ,\R)
\end{equation*}
with 
\begin{equation*}
\tilde{\OO}_{\delta,T_0} := \left\{ \big( \psi^1_0,\psi^2_0,\psi^3_0 \big) \in (\mathcal{S} \cap H^3_{(0)}(0,1))^3 \: ; \: 
\sum_{j=1}^3 ||\psi^j_0 - \Phi_j(T_0)||_{H^3_{(0)}} < \delta \right\},
\end{equation*}
such that $\tilde{\Gamma}_{T_0,T_f} \big( \Phi_1(T_0), \Phi_2(T_0), \Phi_3(T_0) \big)=0$ and, if $(\psi^1_0,\psi^2_0,\psi^3_0) \in \tilde{\OO}_{T_0,\delta}$, the solution $\big( \psi^1, \psi^2, \psi^3 \big)$ of system~(\ref{syst_Neq}) with initial conditions $\psi^j(T_0,\cdot)=\psi^j_0$, for $j =1,2,3$, and control $u:= \tilde{\Gamma}_{T_0,T_f} \big( \psi^1_0,\psi^2_0,\psi^3_0 \big)$ satisfies
\begin{align}
\label{uref_condition3}
&\PP_1 \big( \psi^1(T_f) \big) = \PP_2 \big( \psi^2(T_f) \big) = \PP_3 \big( \psi^3(T_f) \big) = 0,
\\
\label{uref_condition4}
&\Im \big( \lag \psi^1(T_f) ,\Phi_1(T_f) \rag^5 \overline{\lag \psi^2(T_f) , \Phi_2(T_f) \rag^8} 
\lag \psi^3(T_f) , \Phi_3(T_f) \rag^3 \big) \Big) = 0.
\end{align}
\end{prop}

\begin{rmq}
The conditions (\ref{uref_condition3}) and (\ref{uref_condition4}) will be used in the next subsection to prove (\ref{3eq_temps_final}). Equation (\ref{uref_condition4}) will be used to define the global phase $\theta^\eta$.
\end{rmq}

\noindent
\textbf{Proof of Proposition~\ref{traj_ref_inversion_locale} :}
Let us define the following space
\begin{equation*}
X_1 := \left\{ (\phi_1, \phi_2 , \phi_3) \in H^3_{(0)}((0,1),\C)^3 \, ; \, \lag \phi_j , \varphi_k \rag =0, \, \text{ for } \, 1 \leq k \leq j \leq 3 \right\}.
\end{equation*}
We consider the following end-point map
\begin{equation*}
\begin{array}{cccl}
\Theta_{T_0,T_f}: &  L^2((T_0,T_f),\mathbb{R}) \times  H^3_{(0)}(0,1)^3
& \rightarrow & H^3_{(0)}(0,1)^3 \times X_1 \times \R,
\end{array}
\end{equation*}
defined by
\begin{align*}
\Theta_{T_0,T_f} \big( u, \psi^1_0, \psi^2_0, \psi^3_0 \big) := 
\Big( \psi^1_0, \psi^2_0 , \psi^3_0, \, \PP_1 \big( \psi^1(T_f) \big), \PP_2 \big( \psi^2(T_f) \big), \PP_3 \big( \psi^3(T_f) \big),
\\
\Im \big( \, \lag \psi^1(T_f) ,\Phi_1(T_f) \rag^5 \overline{\lag \psi^2(T_f) , \Phi_2(T_f) \rag^8} 
\lag \psi^3(T_f) , \Phi_3(T_f) \rag^3 \, \big) \Big)
\end{align*}
where $(\psi^1,\psi^2,\psi^3)$ is the solution of (\ref{syst_Neq}) with initial condition $\psi^j(T_0,\cdot)=\psi^j_0$ and control $u$. Thus, we have
\begin{equation*}
\Theta_{T_0,T_f} \big( 0, \Phi_1(T_0), \Phi_2(T_0), \Phi_3(T_0) \big) = \big( \Phi_1(T_0), \Phi_2(T_0), \Phi_3(T_0) , 0 , 0 , 0 , 0 \big).
\end{equation*}
Proposition~\ref{traj_ref_inversion_locale} is proved by application of the inverse mapping theorem to $\Theta_{T_0,T_f}$ at the point $\big( 0, \Phi_1(T_0), \Phi_2(T_0), \Phi_3(T_0) \big)$.


\noindent
Using the same arguments as in \cite[Proposition 3]{BeauchardLaurent}, it comes that $\Theta_{T_0,T_f}$ is a $C^1$ map and that
\begin{align*}
\md & \Theta_{T_0,T_f} \big( 0, \Phi_1(T_0), \Phi_2(T_0), \Phi_3(T_0) \big) . (v , \Psi^1_0 , \Psi^2_0 ,\Psi^3_0) 
\\
=& \Big( \Psi^1_0, \Psi^2_0, \Psi^3_0 , \,
\PP_1 \big(\Psi^1(T_f) \big), \PP_2 \big(\Psi^2(T_f) \big), \PP_3 \big(\Psi^3(T_f) \big) ,
\\
& 5 \Im(\lag \Psi^1(T_f),\Phi_1(T_f)\rag) - 8 \Im(\lag \Psi^2(T_f),\Phi_2(T_f)\rag) + 3 \Im(\lag \Psi^3(T_f),\Phi_3(T_f)\rag) \Big),
\end{align*}
where $(\Psi^1, \Psi^2, \Psi^3)$ is the solution of (\ref{linearise_u=0}) on the time interval $(T_0,T_f)$ with control $v$ and initial conditions $\Psi^j(T_0,\cdot) = \Psi^j_0$.

\noindent
It remains to prove that $\md \Theta_{T_0,T_f} \big( 0,\Phi_1(T_0), \Phi_2(T_0) , \Phi_3(T_0) \big):L^2((T_0,T_f),\R) \times H^3_{(0)}(0,1)^3 \to H^3_{(0)}(0,1)^3 \times X_1 \times \R \,$ admits a continuous right inverse. 

\noindent
Let $(\Psi^1_0,\Psi^2_0,\Psi^3_0) \in H^3_{(0)}(0,1)^3$,  $(\psi^1_f,\psi^2_f,\psi^3_f) \in X_1$ and $r \in \R$. 
Straightforward computations lead to
\begin{equation*}
\Psi^j(T_f) =  \sum_{k=1}^{+ \infty} \Big( \lag \Psi^j_0 , \Phi_k(T_0) \rag + i\lag \mu \varphi_j , \varphi_k \rag  
\int_{T_0}^{T_f} v(t) e^{i(\lambda_k - \lambda_j) t} \md t \Big) \Phi_k(T_f).
\end{equation*}
Finding $v \in L^2((T_0,T_f),\R)$ such that
\begin{align*}
&\PP_j(\Psi^j(T_f)) = \psi^j_f, \quad \forall j \in \{1,2,3\},
\\
& \Im \big( 5 \lag \Psi^1(T_f) ,\Phi_1(T_f)\rag - 8 \lag \Psi^2(T_f),\Phi_2(T_f)\rag + 3 \lag \Psi^3(T_f),\Phi_3(T_f)\rag \big) = r,
\end{align*}
is equivalent to solving the following trigonometric moment, $\forall j=1,2,3,  \, \forall k \geq j+1$
\begin{equation}
\label{traj_ref_inversion_locale_pb_moment}
\begin{aligned}
&\int_{T_0}^{T_f} v(t) e^{i (\lambda_k- \lambda_j) t} \md t = \frac{1}{i \lag \mu \varphi_j, \varphi_k \rag} \big( \lag \psi^j_f , \Phi_k(T_f) \rag - \lag \Psi^j_0 , \Phi_k(T_0) \rag \big), 
\\
&\int_{T_0}^{T_f} v(t) \md t = \frac{r - \Im\big( 5 \lag \Psi^1_0,\Phi_1(T_0) \rag - 8 \lag \Psi^2_0,\Phi_2 (T_0) \rag + 3 \lag \Psi^3_0,\Phi_3(T_0) \rag \big)}{5 \lag \mu \varphi_1, \varphi_1 \rag - 8 \lag \mu \varphi_2 ,\varphi_2 \rag + 3 \lag \mu \varphi_3 ,\varphi_3 \rag}.
\end{aligned}
\end{equation}
Using Proposition~\ref{prop_pb_moment} and the hypotheses on $\mu$, this ends the proof of Proposition~\ref{traj_ref_inversion_locale}.
\hfill $\blacksquare$

\subsection{Proof of Theorem~\ref{th_construction_uref}}

\noindent
Let $\delta >0$ be the radius defined in Proposition~\ref{traj_ref_inversion_locale} with $T_0=\varepsilon$ and $T_f=T_1$. 
For $\eta > 0$ we define the following control
\begin{equation}
\label{def_uref}
u_{ref}^{\eta} (t) :=
\left\{
\begin{aligned}
& 0 \quad & \text{for } t \in (0,\frac{\varepsilon}{2}),
\\
\hat{\Gamma} & (\eta) \quad & \text{for } t \in (\frac{\varepsilon}{2},\varepsilon),
\\
\tilde{\Gamma}_{\varepsilon,T_1}(\psi^{1,\eta}_{ref}(\varepsilon)&, \psi^{2,\eta}_{ref}(\varepsilon) , \psi^{3,\eta}_{ref}(\varepsilon)) \quad & \text{for } t \in (\varepsilon,T_1),
\end{aligned}
\right.
\end{equation}
where $\hat{\Gamma}$ and $\tilde{\Gamma}$ are defined respectively in Proposition~\ref{conditions_minimalite} and \ref{traj_ref_inversion_locale}.
We prove that, for $\eta$ small enough, this control satisfies the conditions of Theorem~\ref{th_construction_uref}.
\\

\noindent
\textbf{Proof of Theorem~\ref{th_construction_uref} :}
The proof is decomposed into two parts. First, we prove that there exists $\overline{\eta} > 0$ such that for $\eta \in (0,\overline{\eta})$, $u_{ref}^{\eta}$ is well defined, satisfies $||u_{ref}^{\eta}||_{L^2(0,T_1)} \leq C \eta$ and the conditions (\ref{uref_condition1}), (\ref{uref_condition2}) are satisfied. Then, we prove the existence of $T^{\eta} >0$ and $\theta^\eta \in \R$ such that if $u_{ref}^{\eta}$ is extended by $0$ on $(T_1 , T^{\eta})$, the condition (\ref{3eq_temps_final}) is satisfied.
\\

\textit{First step :}  $u_{ref}^{\eta}$ is well defined.

\noindent
Using Proposition~\ref{conditions_minimalite}, the control $u_{ref}^{\eta}$ is well defined on $(0, \varepsilon)$ as soon as $ \eta \in (0,\eta^*)$.
Moreover, using Lipschitz property of $\hat{\Gamma}$, there exists $C(\eta^*)>0$ such that
\begin{equation*}
|| u_{ref}^{\eta} ||_{L^2(\frac{\varepsilon}{2},\varepsilon)} 
= || \hat{\Gamma} (\eta) - \hat{\Gamma}(0) ||_{L^2(\frac{\varepsilon}{2},\varepsilon)}
\leq C(\eta^*) \eta.
\end{equation*}
Thanks to Proposition~\ref{bien_pose}, there exists $C(\varepsilon)>0$ such that if $||u||_{L^2(0, \varepsilon)} < 1$, the associated solution of (\ref{syst_Neq})-(\ref{cond_ini}) satisfies 
\begin{equation*}
|| ( \psi^j - \Phi_j ) (\varepsilon)||_{H^3_{(0)}} \leq C(\varepsilon) ||u||_{L^2(0,\varepsilon)}, 
\quad \text{for } j=1,2,3.
\end{equation*}
Thus, using Proposition~\ref{traj_ref_inversion_locale}, if $C(\varepsilon) C(\eta^*) \eta < \frac{\delta}{3}$, we get that for $j=1,2,3$,
\begin{equation*}
||  ( \psi^{j,\eta}_{ref} - \Phi_j ) (\varepsilon)||_{H^3_{(0)}} < \frac{\delta}{3}.
\end{equation*}
Thus, $u_{ref}^{\eta}$ is well defined on $(0,T_1)$. Moreover, there exists $C(\delta)>0$ such that
\begin{align*}
&|| u_{ref}^{\eta} ||_{L^2(\varepsilon,T_1)}
\\
&= || \tilde{\Gamma}_{\varepsilon,T_1} \big( \psi^{1,\eta}_{ref}(\varepsilon),\psi^{2,\eta}_{ref}(\varepsilon),\psi^{3,\eta}_{ref}(\varepsilon) \big) - \tilde{\Gamma}_{\varepsilon,T_1} \big( \Phi_1(\varepsilon) ,\Phi_2(\varepsilon) ,\Phi_3(\varepsilon) \big) ||_{L^2(\varepsilon,T_1)}
\\
&\leq C(\delta) \sum_{j=1}^3 || (\psi^{j,\eta}_{ref}-\Phi_j) (\varepsilon)||_{H^3_{(0)}}
\\
&\leq 3 C(\delta) C(\varepsilon) C(\eta^*) \eta.
\end{align*}
Finally, choosing
\begin{equation*}
\overline{\eta} < \min \left( \eta^* , \frac{\delta}{3 C(\varepsilon) C(\eta^*)} , \frac{1}{C(\eta^*)} \right),
\end{equation*}
implies that $||u_{ref}^{\eta}||_{L^2(0,T_1)} \leq C \eta$. Here and throughout this paper $C$ denotes a positive constant that may vary each time it appears. Thanks to Proposition~\ref{conditions_minimalite}, it comes that (\ref{uref_condition1}) and (\ref{uref_condition2}) hold.
\\

\textit{Second step :} We prove the existence of a final time $T^{\eta}>0$ and a global phase $\theta^\eta \in \R$ such that (\ref{3eq_temps_final}) holds.

\noindent
Proposition~\ref{traj_ref_inversion_locale}, implies
\begin{align}
\label{traj_ref_3eq_condition1}
\psi^{j,\eta}_{ref}(T_1) &= \sum_{k=1}^j \lag \psi^{j,\eta}_{ref}(T_1) , \Phi_k(T_1) \rag \Phi_k(T_1) , \quad \forall  j=1,2,3,
\\
\label{traj_ref_3eq_condition2}
\Im \big( \lag \psi^{1,\eta}_{ref}&(T_1) ,\Phi_1(T_1) \rag^5 \overline{\lag \psi^{2,\eta}_{ref}(T_1) , \Phi_2(T_1) \rag^8} \lag \psi^{3,\eta}_{ref}(T_1) , \Phi_3(T_1) \rag^3 \big) =0.
\end{align}
Using the invariant of the system, $\lag \psi^{j,\eta}_{ref} , \psi^{k,\eta}_{ref} \rag \equiv \delta_{j=k}$, for $j,k \in \{1,2,3\}$, this leads to the existence of $\theta^\eta_1$, $\theta^\eta_2, \theta^\eta_3 \in (-\pi ,\pi]$ such that 
\begin{equation*}
\psi^{j,\eta}_{ref}(T_1) = e^{- i \theta^\eta_j} \Phi_j(T_1), \quad \forall j=1,2,3.
\end{equation*}
Using (\ref{traj_ref_3eq_condition2}), it comes that
\begin{equation*}
\sin \big(5 \theta^\eta_1 - 8 \theta^\eta_2 + 3 \theta^\eta_3 \big) = 0.
\end{equation*}
Using Proposition~\ref{bien_pose}, it comes that, up to a choice of a smaller $\overline{\eta}$,
\begin{equation}
\label{traj_ref_3eq_condition_phase}
5 \theta^\eta_1 - 8 \theta^\eta_2 + 3 \theta^\eta_3 = 0.
\end{equation}
Recall that $\lambda_k = k^2 \pi^2$. Let $T^{\eta}$ and $\theta^\eta$ be such that $T^{\eta}>T_1$ and
\begin{equation*}
\left\{ 
\begin{aligned}
&T^{\eta}  \equiv \frac{\theta^\eta_1 - \theta^\eta_2}{\lambda_2-\lambda_1} \:  \left[ \frac{2}{\pi} \right],
\\
&\theta^\eta \equiv \frac{\lambda_2}{\lambda_2 - \lambda_1} \theta^\eta_1 - \frac{\lambda_1}{\lambda_2-\lambda_1} \theta^\eta_2 \: 
[2 \pi].
\end{aligned}
\right.
\end{equation*}
This choice leads to
\begin{equation*}
\left\{
\begin{aligned}
\theta^\eta_1 + \lambda_1 T^{\eta} - \theta^\eta & \equiv 0 \: [2 \pi],
\\
\theta^\eta_2 + \lambda_2 T^{\eta} - \theta^\eta & \equiv 0 \: [2 \pi].
\end{aligned}
\right.
\end{equation*}
Then, using the definitions of $T^{\eta}$ and $\theta^\eta$ together with (\ref{traj_ref_3eq_condition_phase}) we get
\begin{align*}
\theta^\eta_3 + \lambda_3 T^{\eta} - \theta^\eta \, &\equiv  \,
\theta^\eta_3 + \frac{\lambda_3}{\lambda_2-\lambda_1}( \theta^\eta_1 - \theta^\eta_2) - \frac{\lambda_2}{\lambda_2-\lambda_1} \theta^\eta_1 + \frac{\lambda_1}{\lambda_2-\lambda_1} \theta^\eta_2 \: [2 \pi]
\\
& \equiv \, \frac{1}{3} \big( 5 \theta^\eta_1 - 8 \theta^\eta_2 + 3 \theta^\eta_3 \big) \:[2\pi]
\\
& \equiv \, 0 \: [2 \pi].
\end{align*}
Finally, if we extend $u_{ref}^{\eta}$ by $0$ on $(T_1,T^{\eta})$, we have that $\big( \psi^{1,\eta}_{ref}, \psi^{2,\eta}_{ref}, \psi^{3,\eta}_{ref} \big)$ is solution of (\ref{syst_Neq})-(\ref{cond_ini}) with control $u_{ref}^{\eta}$ and satisfies for $j \in \{1,2,3\}$
\begin{equation*}
\psi^{j,\eta}_{ref}(T^{\eta}) = e^{-i( \theta^\eta_j + \lambda_j T^{\eta} )} \varphi_j = e^{-i \theta^\eta} \varphi_j.
\end{equation*}
This ends the proof of Theorem~\ref{th_construction_uref}.

\hfill $\blacksquare$

\section{Proof of Theorem~\ref{theo_controle_3eq}}
\label{sect_methode_retour}

\noindent
This section is dedicated to the proof of Theorem~\ref{theo_controle_3eq} which is done in the case $T=0$, the extension to the general case being straightforward. The proof is divided in two parts. In Subsection \ref{subsect_controle_L}, the functional setting is specified and we prove the controllability of the linearized system around $(\psi^{1,\eta}_{ref}, \psi^{2,\eta}_{ref}, \psi^{3,\eta}_{ref}, u_{ref}^{\eta})$,
\begin{equation}
\label{linearise_ref}
\left\{
\begin{aligned}
&i \partial_t \Psi^{j,\eta} = -\partial^2_{xx} \Psi^{j,\eta} - u_{ref}^{\eta}(t) \mu(x) \Psi^{j,\eta} - v(t) \mu(x) \psi^{j,\eta}_{ref}, 
\: &(t,x)& \in (0,T^{\eta}) \times (0,1), 
\\
&\Psi^{j,\eta}(t,0) = \Psi^{j,\eta}(t,1) =0 , \quad &t& \in (0,T^{\eta}), 
\\
&\Psi^{j,\eta}(0,x) = 0, \quad &x& \in (0,1),
\end{aligned}
\right.
\end{equation}
when $\eta$ is small enough.
In Subsection \ref{subsect_controle_NL}, we conclude the proof of Theorem~\ref{theo_controle_3eq} using the inverse mapping theorem.

\subsection{Controllability of the linearized system}
\label{subsect_controle_L}

\noindent
For any $t>0$, let
\begin{equation}
\label{def_XfT}
\begin{aligned}
X^f_t :&= \Big\{ (\phi^1,\phi^2,\phi^3) \in H^3_{(0)}((0,1),\C)^3 \, ; \,
\Re (\lag \phi^j , \psi^{j,\eta}_{ref}(t) \rag)=0, \text{ for } j=1,2,3  \\
&\text{ and } \lag \phi^j , \psi^{k,\eta}_{ref}(t) \rag = - \overline{\lag \phi^k , \psi^{j,\eta}_{ref}(t)\rag}, \text{ for } 
(j,k)=(2,1),(3,1),(3,2) \Big\}.
\end{aligned}
\end{equation}
The following proposition holds. 
\begin{prop}
\label{controle_linearise}
There exists $\hat{\eta} \in (0,\overline{\eta})$ such that, for any $\eta \in (0,\hat{\eta})$, if $T^{\eta}$, $u_{ref}^{\eta}$ and $(\psi^{1,\eta}_{ref}, \psi^{2,\eta}_{ref}, \psi^{3,\eta}_{ref})$ are defined as in Theorem~\ref{th_construction_uref}, there exists a continuous linear map
\begin{equation*}
\begin{array}{cccc}
L^{\eta}: &  X^f_{T^{\eta}}
& \rightarrow & L^2((0,T^{\eta}),\mathbb{R})
\\
& (\psi^1_f,\psi^2_f,\psi^3_f) & \mapsto & v
\end{array}
\end{equation*}
such that for any $(\psi^1_f,\psi^2_f,\psi^3_f) \in X^f_{T^{\eta}}$, the solution $(\Psi^{1,\eta}, \Psi^{2,\eta}, \Psi^{3,\eta})$ of system~(\ref{linearise_ref}) with control $v=L^{\eta}(\psi^1_f,\psi^2_f,\psi^3_f)$ satisfies
\begin{equation*}
\big( \Psi^{1,\eta}(T^{\eta}), \Psi^{2,\eta}(T^{\eta}), \Psi^{3,\eta}(T^{\eta}) \big) = \big( \psi^1_f, \psi^2_f, \psi^3_f \big).
\end{equation*}
\end{prop}

\noindent
Before proving Proposition~\ref{controle_linearise} we set some notations. 
For any $\eta \in (0,\overline{\eta})$, for any $t \in (0,T^{\eta})$, let $U^\eta(t)$ be the propagator of the following system
\begin{equation} \label{syst_1eq_ref}
\left\{ 
\begin{aligned}
& i \partial_t \psi = - \partial^2_{xx} \psi - u_{ref}^{\eta} (t) \mu(x) \psi, \: &(t,x)& \in (0,T^{\eta})\times(0,1), 
\\
& \psi (t,0) = \psi(t,1) = 0, \: &t& \in (0,T^{\eta}),
\\
& \psi (0,x) = \psi^0(0,x), \: &x& \in (0,1),
\end{aligned}
\right.
\end{equation}
i.e. $U^{\eta}(t) \psi^0 = \psi(t)$. We will work in the Hilbert basis $(\Phi^{\eta}_k(t) := U^{\eta}(t) \varphi_k)_{k \in \N^*}$ of $L^2((0,1),\C)$. Notice that for $j=1,2,3$, $\Phi^{\eta}_j = \psi^{j,\eta}_{ref}$.
As the proof of Proposition~\ref{controle_linearise} is quite long and technical, let us detail the different steps. Let
\begin{equation*}
\II := \big\{ (j,k) \in \{1,2,3 \} \times \N^* \, ; \, k \geq j+1 \big\} \cup \left\{ (3,3) \right\}.
\end{equation*}
\begin{sloppypar} 
The first step consists in proving the controllability of the components $\lag \Psi^{j,\eta}(T_f),\Phi^{\eta}_k(T_f) \rag$ for $(j,k) \in \II$, for any $T_f >0$ and $\eta$ sufficiently small, as stated in Lemma~\ref{lemme_surj2}. First, we prove that these components are controllable when $\eta =0$ : it corresponds to solving a trigonometric moment problem with an infinite asymptotic gap between successive frequencies. Then, we extend the controllability of these components to small values of $\eta$, by an argument of close linear maps.

In the second step (Lemmas \ref{lemme_base_riesz} and \ref{lemme_famille_minimale}), using Riesz basis and biorthogonal family arguments, we prove that we can also control the two diagonal directions $\lag \Psi^{j,\eta}(T_f),\Phi^{\eta}_j(T_f) \rag$ for $j=1,2$. This would not have been possible directly in the first step. Indeed for $\eta=0$, the three directions $\lag \Psi^{j,\eta}(T_f), \Phi^{\eta}_j(T_f) \rag$ for $j=1,2,3$ are associated to the same frequency in the moment problem. But for $\eta >0$, the construction of the reference trajectory (and more precisely conditions (\ref{uref_condition1}) and (\ref{uref_condition2})) will allow to control those two directions.

Finally, in the third step, due to the conditions imposed in the definition of $X^f_t$ (in (\ref{def_XfT})) the remaining directions $\lag \Psi^{j,\eta} , \Phi^{\eta}_k \rag$ for $1 \leq k < j$ are automatically controlled.
\\
\end{sloppypar}

\noindent
\textbf{Proof of Proposition~\ref{controle_linearise} :}

The map $L^\eta$ will be designed on $(0,T_1)$ and extended by $0$ on $(T_1,T^{\eta})$, where $T_1$ is as in Theorem~\ref{th_construction_uref}. Let
\begin{equation*}
\VV_0 := \left\{ (d^1,d^2,d^3) \in h^3(\N^*,\C)^3 \, ; \, d^j_k =0 , \text{ if } (j,k) \notin \II \text{ and } \Re(d^3_3)=0 \right\}.
\end{equation*}
Let $R : \II \to \N$ be the rearrangement such that, if $\omega_n := \lambda_k- \lambda_j$ with $n=R(j,k)$, the sequence $(\omega_n)_{n \in \N}$ is increasing. Notice that $0=R(3,3)$.
\\

\textit{First step of the proof of Proposition~\ref{controle_linearise} :} we prove that the directions $\lag \Psi^{j,\eta}(T_f),\Phi^{\eta}_k(T_f) \rag$ for $(j,k) \in \II$ are controllable in any positive time $T_f$ for $\eta$ small enough.

\noindent
Let
\begin{equation*}
d^\eta_{T_f} : \psi=(\psi^1,\psi^2,\psi^3) \in X^f_{T_f} \mapsto \big(d_{T_f}^{1,\eta}(\psi), d_{T_f}^{2,\eta}(\psi), d_{T_f}^{3,\eta}(\psi) \big) \in \VV_0,
\end{equation*}
where for $j=1,2,3$, 
\begin{align*}
d^{j,\eta}_{T_f ,k}(\psi) :&= \lag \psi^j , \Phi^{\eta}_k(T_f) \rag, \quad \text{if } (j,k) \in \II,
\\
d^{j,\eta}_{T_f ,k}(\psi) :&= 0, \quad \text{if } (j,k) \not\in \II.
\end{align*}
The next lemma ensures the controllability of the directions $\lag \Psi^{j,\eta}(T_f),\Phi^{\eta}_k(T_f) \rag$ for $(j,k) \in \II$.
\begin{lemme}
\label{lemme_surj2}
Let $T_f > 0$ and
\begin{equation*}
\begin{array}{cccc}
F^\eta : & L^2((0,T_f),\mathbb{R}) 
& \rightarrow &  \VV_0
\\
& v & \mapsto & d^\eta_{T_f}(\Psi(T_f))
\end{array}
\end{equation*}
where $\Psi := \big( \Psi^1, \Psi^2, \Psi^3 \big)$ is the solution of (\ref{linearise_ref}) with control $v$.
There exists $\hat{\eta}=\hat{\eta}(T_f) \in (0,\overline{\eta})$ such that, for any $\eta \in (0,\hat{\eta})$, the map $F^\eta$ has a continuous right inverse 
\begin{equation*}
F^{\eta^{-1}} : \VV_0 \rightarrow L^2((0,T_f),\R).
\end{equation*}
\end{lemme}

\noindent
\textit{Proof of Lemma~\ref{lemme_surj2} :}
Straightforward computations lead to
\begin{equation} \label{expression_linearise}
\lag \Psi^{j,\eta}(T_f) , \Phi^{\eta}_k(T_f) \rag
= i \int_0^{T_f} v(t)\lag \mu \psi^{j,\eta}_{ref}(t), \Phi^{\eta}_k(t) \rag \md t, \quad \text{for } (j,k) \in \II.
\end{equation}
Let us define
\begin{equation} \label{def_fn}
f^{\eta}_n(t) := \frac{\lag \mu \psi^{j,\eta}_{ref}(t) , \Phi^{\eta}_k(t) \rag}{\lag \mu \varphi_j , \varphi_k \rag}, \; \text{ for } (j,k) \in \II \text{ and } n=R(j,k),
\end{equation}
and $f^{\eta}_{-n}(t) := \overline{f^{\eta}_n(t)}$, for $n \in \N^*$. We consider the following map
\begin{equation*}
\begin{array}{cccc}
J^{\eta}: &  L^2((0,T_f),\C)
& \rightarrow & \ell^2(\Z,\C)
\\
& v & \mapsto &  \left( \int_0^{T_f} v(t) f^{\eta}_n(t) \md t \right)_{n \in \Z}
\end{array}.
\end{equation*}
Notice that $f_n^0(t)=e^{i \omega_n t}$ with $\omega_n = \lambda_k - \lambda_j$ for any $n=R(j,k) \in \N$.
Thus (see \ref{annexe_pb_moment}), $J^0$ is continuous with values in $\ell^2(\Z,\C)$. Moreover, $J^0$ is an isomorphism from $H_0 :=  \text{Adh}_{L^2(0,T_f)} \left( \text{Span} \{ f_n^0 \, ;\, n \in \Z\} \right)$ to $\ell^2(\Z,\C)$.
\\

\textit{First step : } we prove the existence of $\tilde{C}>0$ such that
\begin{equation} \label{applications_proches}
|| (J^\eta - J^0)(v) ||_{\ell^2} \leq \tilde{C} \eta || v ||_{L^2(0,T_f)}, \quad \forall v \in L^2((0,T_f),\C).
\end{equation}

\noindent
Let $(j,k) \in \II$ and $n=R(j,k)\in \N$ and $v \in L^2((0,T_f),\C)$. Using (\ref{expression_linearise}) and (\ref{def_fn}), the triangular inequality and Hypothesis~\ref{hypo_mu}, we get
\begin{align*}
&\left| \int_0^{T_f} v(t) (f_n^0 - f_n^{\eta})(t) \md t \right|
= \left|  \frac{\langle \Psi^{j,0}(T_f), \Phi_k(T_f) \rangle}{\lag \mu \varphi_j ,\varphi_k \rag} 
-  \frac{\langle \Psi^{j,\eta}(T_f), \Phi^{\eta}_k(T_f) \rangle}{\lag \mu \varphi_j ,\varphi_k \rag} \right|
\\
&\leq C k^3 \left(\left| \langle (\Psi^{j,0} - \Psi^{j,\eta})(T_f) , \Phi_k(T_f) \rangle \right|
+ \left| \langle (U^{\eta}(T_f) - U^0(T_f))^* \Psi^{j,\eta}(T_f), \varphi_k \rangle \right| \right)
\end{align*}
because $(\Phi^{\eta}_k - \Phi_k)(t) = (U^{\eta}(t) - U^0(t)) \varphi_k$ (we denoted by $*$ the $L^2((0,1),\C)$ adjoint operator).
Thus,
\begin{align*} 
|| (J^0-J^{\eta})(v) ||_{\ell^2} \leq C
&\sum_{j=1}^3 \Big( || (\Psi^{j,0} - \Psi^{j,\eta}) (T_f) ||_{H^3_{(0)}} +
\\
&  || (U^{\eta}(T_f) - U^0(T_f))^* \Psi^{j,\eta}(T_f) ||_{H^3_{(0)}} \Big).
\label{applications_proches1}
\tag{\theequation} \addtocounter{equation}{1}
\end{align*}
Proposition~\ref{bien_pose} implies that
\begin{align*}
&|| (\Psi^{j,0} - \Psi^{j,\eta})(T_f)||_{H^3_{(0)}} \\
&\leq C ||u_{ref}^{\eta}(t) \mu \Psi^{j,0}(t) + v(t) \mu (\psi^{j,\eta}_{ref}-\Phi_j)(t)||_{L^2((0,T_f), H^3\cap H^1_0)}
\\
&\leq C ||u_{ref}^{\eta} ||_{L^2(0,T_f)} ||v||_{L^2(0,T_f)}.
\label{applications_proches2}
\tag{\theequation} \addtocounter{equation}{1}
\end{align*}

\noindent
Using unitarity, it comes that $U^{\eta}(T_f)^*$ is the propagator at time $T_f$ of system
\begin{equation*}
\left\{
\begin{aligned}
&i \partial_t \psi = \partial^2_{xx} \psi + u_{ref}^{\eta}(T_f -t) \mu(x) \psi, \: &(t,x)& \in (0,T_f) \times (0,1),
\\
& \psi(t,0) = \psi(t,1) =0, \: &t& \in (0,T_f).
\end{aligned}
\right.
\end{equation*}
Thus Proposition~\ref{bien_pose} may be applied again leading to
\begin{align*}
&\big| \big| (U^{\eta}(T_f)- U^0(T_f))^* \Psi^{j,\eta}(T_f) \big| \big|_{H^3_{(0)}}
\\
&\leq C \, || u_{ref}^{\eta} (t) \mu U^0(t)^* \Psi^{j,\eta}(T_f) ||_{L^2((0,T_f),H^3 \cap H^1_0)}
\\
&\leq C \, ||u_{ref}^{\eta} ||_{L^2(0,T_f)} || v ||_{L^2(0,T_f)}.
\label{applications_proches3}
\tag{\theequation} \addtocounter{equation}{1}
\end{align*}
From inequalities (\ref{applications_proches1}), (\ref{applications_proches2}), (\ref{applications_proches3}) above and (\ref{borne_uref}) we get the conclusion of the first step.
\\

\textit{Second step : } conclusion.

\noindent
Let $\hat{\eta}(T_f) := \min \left\{ \overline{\eta}, \, \tilde{C}^{-1} || (J^0)^{-1} ||_{\LL(H_0,\ell^2)}^{-1} \right\}$ where $\tilde{C}$ is defined by (\ref{applications_proches}) and let $\eta \in (0,\hat{\eta}(T_f))$.
We deduce from the first step that $J^{\eta}$ is an isomorphism from $H_0$ to $\ell^2(Z,\C)$. Let $(d^1,d^2,d^3) \in \VV_0$. We define
$\tilde{d}_n := \dfrac{d^j_k}{i\lag \mu \varphi_j, \varphi_k \rag}$, for $(j,k) \in \II$ and $n = R(j,k) \in \N$,
and $\tilde{d}_{-n} := \overline{\tilde{d}_n}$, for $n \in \N^*$. Then, 
\begin{equation*}
F^{\eta^{-1}}(d^1,d^2,d^3) := (J^{\eta}_{|H_0})^{-1} ( \tilde{d})
\end{equation*}
is the unique solution $v$ in $H_0$ of the equation $F^\eta(v) = (d^1,d^2,d^3)$. The uniqueness implies that $v$ is real valued.
This ends the proof of Lemma~\ref{lemme_surj2}. 

\hfill $\blacksquare$ \\

\textit{Second step of the proof of Proposition~\ref{controle_linearise} :} Riesz basis and minimality.

\noindent
To prove that we can also control the directions $\lag \Psi^{j,\eta}(T_f) , \Phi^{\eta}_j(T_f) \rag$, for $j=1,2$, we will use the following lemmas.
\begin{lemme}
\label{lemme_base_riesz}
Let $T_f >0$ and $H^{\eta} := \text{Adh}_{L^2(0,T_f)}\big( \text{Span} \{f^{\eta}_n , n \in \Z\} \big)$. If $\eta < \hat{\eta}(T_f)$, then $(f^{\eta}_n)_{n \in \Z}$ is a Riesz basis of $H^{\eta}$.
\end{lemme}

\noindent
\textit{Proof of Lemma~\ref{lemme_base_riesz} : }
Using \cite[Proposition 19]{BeauchardLaurent}, it comes that $(f^{\eta}_n)_{n \in \Z}$ is a Riesz basis of $H^{\eta}$ if and only if there exists $C_1, C_2 >0$ such that for any complex sequence $(a_n)_{n \in \Z}$ with finite support
\begin{equation} \label{base_riesz}
C_1 \left( \sum_n |a_n|^2 \right)^{1/2} 
\leq 
\left| \left|  \sum_n a_n f^{\eta}_n \right| \right|_{L^2(0,T_f)}
\leq
C_2 \left( \sum_n |a_n|^2 \right)^{1/2} .
\end{equation}
Lemma~\ref{lemme_surj2} together with \cite[Theorem 1]{Boas} imply the first inequality of (\ref{base_riesz}).
Using again \cite[Theorem 1]{Boas}, we get that the second inequality of (\ref{base_riesz}) holds 
if and only if, for any $g \in L^2((0,T_f),\C)$
\begin{equation*}
\left( \sum_{n \in \Z} \Big| \int_0^{T_f} g(t)  f^{\eta}_n(t) \md t \Big|^2 \right)^{1/2} \leq C_2 ||g||_{L^2}.
\end{equation*}
This is implied by the continuity of $J^0$, the triangular inequality and (\ref{applications_proches}). This ends the proof of Lemma~\ref{lemme_base_riesz}.

\hfill $\blacksquare$
\\

\noindent
From now on, we consider $\hat{\eta} < \min \big( \hat{\eta}(\frac{\varepsilon}{2}), \hat{\eta}(T_1) \big)$ and $\eta \in (0,\hat{\eta})$ fixed for all what follows.
\\

\begin{lemme}
\label{lemme_famille_minimale}
Let $f^{\eta}_{j,j} := \dfrac{\lag \mu \psi^{j,\eta}_{ref} , \psi^{j,\eta}_{ref} \rag}{\lag \mu \varphi_j, \varphi_j \rag}$, for $j\in \{1,2\}$. The family $\Xi := \big( f^{\eta}_n \big)_{n \in \Z} \cup \{ f^{\eta}_{1,1} , f^{\eta}_{2,2} \}$ is minimal in $L^2((0,T_1),\C)$. 
\end{lemme}

\noindent
\textit{Proof of Lemma~\ref{lemme_famille_minimale} : }
Assume that there exist $(c_n)_{n \in \Z} \in \ell^2(\Z,\C)$ and $c_{1,1}$, $c_{2,2} \in \C$, not all being zero, such that
\begin{equation}
\label{famille_minimale1}
c_{1,1} f^{\eta}_{1,1} + c_{2,2} f^{\eta}_{2,2} + \sum_{n \in \Z} c_n f^{\eta}_n= 0 , \quad \text{in } L^2((0,T_1),\C).
\end{equation}
Thus,
\begin{equation*}
c_{1,1} f^{\eta}_{1,1} + c_{2,2} f^{\eta}_{2,2} + \sum_{n \in \Z} c_n f^{\eta}_n= 0 , \quad \text{in } L^2((0,\frac{\varepsilon}{2}),\C).
\end{equation*}
As $f^{\eta}_0 = f^{\eta}_{1,1} = f^{\eta}_{2,2} = 1$ on $(0, \frac{\varepsilon}{2})$, then
\begin{equation*}
c_{1,1} f^{\eta}_{1,1} + c_{2,2} f^{\eta}_{2,2} + c_0 f^{\eta}_0 = c f^{\eta}_0, \quad \text{on } (0,\frac{\varepsilon}{2}),
\end{equation*}
where $c := c_{1,1} + c_{2,2} + c_0$. Thus,
\begin{equation*}
c f^{\eta}_0 + \sum_{n \in \Z^*} c_n f^{\eta}_n =0, \quad \text{in } L^2((0,\frac{\varepsilon}{2}),\C).
\end{equation*}
As $\eta < \hat{\eta}(\varepsilon/2)$, Lemma~\ref{lemme_base_riesz} with $T_f = \varepsilon/2$ implies minimality of $(f^{\eta}_n)_{n \in \Z}$ in $L^2((0,\frac{\varepsilon}{2}),\C)$. Thus,
\begin{equation*}
c=0 \quad \text{and} \quad c_n = 0, \quad \forall n \in \Z^*.
\end{equation*}
Then, equation (\ref{famille_minimale1}) implies that,
\begin{equation}
\label{famille_minimale2}
c_{1,1} f^{\eta}_{1,1} + c_{2,2} f^{\eta}_{2,2} + c_0 f^{\eta}_0 =0 , \text{ on } (0,T_1).
\end{equation}
Finally, as $c=0$, conditions (\ref{uref_condition1}) and (\ref{uref_condition2}) in (\ref{famille_minimale2}) lead to $c_{1,1} = c_{2,2} = 0$ and then $c_0=0$. 
This is a contradiction, thus the family $\Xi$ is proved to be minimal in $L^2((0,T_1),\C)$.

\hfill $\blacksquare$
\\

\noindent
The proof of Lemma~\ref{lemme_famille_minimale} makes important use of the conditions (\ref{uref_condition1}) and (\ref{uref_condition2}) from the construction of the reference trajectory. This is the main interest of the construction of the reference trajectory : for $\eta =0$, one gets $f^{0}_{1,1} = f^{0}_{2,2} = f^{0}_0$. Thus, one could not control simultaneously $\lag \Psi^{j,0}(T_1), \Phi_j(T_1) \rag$ for $j=1,2,3$. In our setting, the minimal family property allows together with Lemma~\ref{lemme_surj2} to conclude the proof of Proposition~\ref{controle_linearise}.
\\

\textit{Third step of the proof of Proposition~\ref{controle_linearise} :} conclusion.

\begin{sloppypar}
\noindent
Using \cite[Proposition 18]{BeauchardLaurent}, Lemma~\ref{lemme_famille_minimale} implies that there exists a unique biorthogonal family associated to $\Xi$ in $\text{Adh}_{L^2(0,T_1)} \big( \text{Span} (\Xi) \big)$  denoted by $\left\{ g^{\eta}_{1,1}, g^{\eta}_{2,2}, (g^{\eta}_n)_{n \in\Z} \right\}$. This construction ensures that $g^{\eta}_{1,1}$ and $g^{\eta}_{2,2}$ are real valued.

\noindent
Let $\psi_f \in X^f_{T^{\eta}}$ and $\tilde{\psi}_f := \big( e^{i A ( T^{\eta}-T_1)} \psi_f^1, \, e^{i A ( T^{\eta}-T_1)} \psi_f^2, \, e^{i A ( T^{\eta}-T_1)} \psi_f^3 \big)$. As $u_{ref}^{\eta}$ is identically equal to $0$ on $(T_1,T^{\eta})$, it comes that $\tilde{\psi}_f \in X^f_{T_1}$. The map $L^{\eta}$ is defined by
\begin{equation*}
L^{\eta} : \psi_f \in  X^f_{T^{\eta}} \mapsto v \in  L^2((0,T^{\eta}),\mathbb{R}),
\end{equation*}
where $v$ is defined on $(0,T_1)$ by
\begin{equation*}
v := v_0 + \sum_{j=1}^2 \left( \Im \big( \lag \tilde{\psi}^j_f, \psi^{j,\eta}_{ref}(T_1) \rag \big) - \int_0^{T_1} v_0(t) \lag \mu \psi^{j,\eta}_{ref}(t) , \psi^{j,\eta}_{ref}(t) \rag \md t \right) g^{\eta}_{j,j}(t),
\end{equation*}
with $v_0 := F^{\eta^{-1}} \big( d_{T_1}( \tilde{\psi}_f) \big)$ and extended by $0$ on $(T_1, T^{\eta})$.
Notice that $L^{\eta}$ is linear and continuous and that as $v_0$, $g^{\eta}_{1,1}$ and $g^{\eta}_{2,2}$ are real valued so is $v$.
\end{sloppypar}

\noindent
Let $(\Psi^1,\Psi^2,\Psi^3)$ be the solution of (\ref{linearise_ref}) with control $v$.
Using the biorthogonal properties, the definition of $v_0$ and Lemma~\ref{lemme_surj2} with $T_f = T_1$ we get that
\begin{equation*}
\lag \Psi^j(T_1) , \Phi^{\eta}_k(T_1) \rag = \lag \tilde{\psi}^j_f , \Phi^{\eta}_k(T_1) \rag, \quad \forall (j,k) \in \II \cup \{ (1,1),(2,2) \}.
\end{equation*}
We check that $v$ also controls the remaining extra-diagonal terms. Straightforward computations give
\begin{equation*}
\lag \Psi^2(T_1) , \Phi^{\eta}_1(T_1) \rag = - \overline{\lag \Psi^1(T_1) , \Phi^{\eta}_2(T_1) \rag}.
\end{equation*}
Yet, by definition of $v$ and $X^f_{T_1}$,
\begin{equation*}
\lag \Psi^1(T_1) , \psi^{2,\eta}_{ref}(T_1) \rag 
= \lag \tilde{\psi}^1_f , \Phi^{\eta}_2(T_1) \rag
= - \overline{ \lag \tilde{\psi}^2_f , \Phi^{\eta}_1(T_1) \rag}.
\end{equation*}
This leads to
\begin{equation*}
\lag \Psi^2(T_1) , \Phi^{\eta}_1(T_1) \rag = \lag \tilde{\psi}^2_f , \Phi^{\eta}_1(T_1) \rag.
\end{equation*}
The same computations hold for $\lag \Psi^3(T_1), \Phi^{\eta}_1(T_1) \rag$ and $\lag \Psi^3(T_1) , \Phi^{\eta}_2(T_1) \rag$. Thus, as $(\Phi^{\eta}_k(T_1))_{k \in \N^*}$ is a Hilbert basis of $L^2((0,T_1),\C)$, it comes that 
\begin{equation*}
\big( \Psi^1(T_1), \Psi^2(T_1), \Psi^3(T_1) \big) = \big( \tilde{\psi}_f^1, \tilde{\psi}_f^2, \tilde{\psi}_f^3 \big).
\end{equation*}
As $v$ is set to zero on $(T_1,T^{\eta})$, this ends the proof of Proposition~\ref{controle_linearise}.
\hfill $\blacksquare$
\\

\subsection{Controllability of the nonlinear system}
\label{subsect_controle_NL}

\noindent
In this subsection, we end the proof Theorem~\ref{theo_controle_3eq}. First, using the inverse mapping theorem and Proposition~\ref{controle_linearise}, we prove in Proposition~\ref{inversion_locale_methode_retour} that we can control the projections associated to the space $X^f_{T^{\eta}}$ (see below for precise statements and notations). Then, using the invariants of the system (see Remark~\ref{rmq_invariant}) we prove that it is sufficient to control those projections.

\noindent
We define
\begin{equation*}
\begin{array}{cccl}
\Lambda : &  L^2((0,T^{\eta}),\R)
& \rightarrow &  X^f_{T^{\eta}} 
\\
& u & \mapsto & \big( \tilde{\PP}_j(\psi^j(T^{\eta}))_{j=1,2,3} \big)
\end{array}
\end{equation*}
where $(\psi^1,\psi^2,\psi^3)$ is the solution of (\ref{syst_Neq})-(\ref{cond_ini}) with control $u$ and $\tilde{\PP}$ is defined by
\begin{align*}
\tilde{\PP}_j(\phi^j) :&= \phi^j - \Re \big( \lag \phi^j, \psi^{j,\eta}_{ref}(T^{\eta}) \rag \big) \, \psi^{j,\eta}_{ref}(T^{\eta}) \,
\\
- &\, \sum_{k=1}^{j-1} \big( \lag \phi^j , \psi^{k,\eta}_{ref}(T^{\eta}) \rag + \lag \psi^{j,\eta}_{ref}(T^{\eta}),\phi^k \rag \big) \, \psi^{k,\eta}_{ref}(T^{\eta}).
\end{align*}
Thanks to this definition, $\Lambda$ takes value in $X^f_{T^{\eta}}$ (defined in (\ref{def_XfT})) and $\Lambda(u^{\eta}_{ref})=(0,0,0)$. As announced, we prove that we can control the projections $\tilde{\PP}_j$. More precisely, we prove the following proposition.
\begin{prop}
\label{inversion_locale_methode_retour}
There exists $\delta>0$ and a $C^1$-map 
\begin{equation*}
\Upsilon : \Omega_{\delta} \to L^2((0,T^{\eta}),\R),
\end{equation*}
with
\begin{equation*}
\Omega_{\delta} := \left\{ \big( \tilde{\psi}^1_f, \tilde{\psi}^2_f , \tilde{\psi}^3_f \big) \in X^f_{T^{\eta}}  \; ; \;  \sum_{j=1}^3 || \tilde{\psi}^j_f ||_{H^3_{(0)}} < \delta \right\}
\end{equation*}
such that $\Upsilon\big( 0, 0 ,0) = u^{\eta}_{ref}$ and for any $\big( \tilde{\psi}^1_f, \tilde{\psi}^2_f , \tilde{\psi}^3_f \big) \in \Omega_{\delta}$, the solution of system~(\ref{syst_Neq})-(\ref{cond_ini}) with control $u:= \Upsilon\big( \tilde{\psi}^1_f, \tilde{\psi}^2_f , \tilde{\psi}^3_f \big)$ satisfies
\begin{equation*}
\big( \tilde{\PP}_1(\psi^1(T^{\eta})) , \tilde{\PP}_2(\psi^2(T^{\eta})) , \tilde{\PP}_3(\psi^3(T^{\eta})) \big) = \big( \tilde{\psi}^1_f, \tilde{\psi}^2_f , \tilde{\psi}^3_f \big).
\end{equation*}
\end{prop}

\noindent
\textbf{Proof of Proposition~\ref{inversion_locale_methode_retour} :}
This proposition is proved by application of the inverse mapping theorem to $\Lambda$ at the point $u^{\eta}_{ref}$.
Using the same arguments as in \cite[Proposition 3]{BeauchardLaurent}, it comes that $\Lambda$ is $C^1$ and for any $v \in L^2((0,T^{\eta}),\R)$,
\begin{equation*}
\md \Lambda(u^{\eta}_{ref}). v = \big( \tilde{\PP}_1(\Psi^1(T^{\eta})) , \tilde{\PP}_2(\Psi^2(T^{\eta})) , \tilde{\PP}_3(\Psi^3(T^{\eta})) \big),
\end{equation*}
where $(\Psi^{j})_{j=1,2,3}$ is the solution of system~(\ref{linearise_ref}) with control $v$. Straightforward computations lead to $\tilde{\PP}_j(\Psi^{j}(T^{\eta}))= \Psi^{j}(T^{\eta})$ and thus
\begin{equation*}
\md \Lambda(u^{\eta}_{ref}). v = \big( \Psi^1(T^{\eta}) , \Psi^2(T^{\eta}) , \Psi^3(T^{\eta}) \big).
\end{equation*}
Proposition~\ref{controle_linearise} proves that $\md \Lambda(u^{\eta}_{ref}) : L^2((0,T^{\eta}),\R) \to X^f_{T^{\eta}}$ admits a continuous right inverse. This ends the proof of Proposition~\ref{inversion_locale_methode_retour}.

\hfill $\blacksquare$

\noindent
\textbf{Proof of Theorem~\ref{theo_controle_3eq} :}
Let $\tilde{\varepsilon} > 0$ and $\big( \psi^1_f , \psi^2_f ,\psi^3_f \big) \in H^3_{(0)}((0,1),\C)^3$ be such that
\begin{equation*}
\lag \psi^j_f , \psi^k_f \rag = \delta_{j=k}
\, \text{ and } \,
\sum_{j=1}^3 || \psi^j_f - \psi^{j,\eta}_{ref}(T^{\eta}) ||_{H^3_{(0)}} < \tilde{\varepsilon}.
\end{equation*}
Let
\begin{equation*}
\big( \tilde{\psi}^1_f, \tilde{\psi}^2_f , \tilde{\psi}^3_f \big) := \big( \tilde{\PP}_1(\psi^1_f) , \tilde{\PP}_2(\psi^2_f) , \tilde{\PP}_3(\psi^3_f) \big).
\end{equation*}
Let $\delta$ be the radius defined in Proposition~\ref{inversion_locale_methode_retour}. 
There exists $\varepsilon_0 > 0$ such that for any $\tilde{\varepsilon} \in (0, \varepsilon_0)$, $\big( \tilde{\psi}^1_f, \tilde{\psi}^2_f , \tilde{\psi}^3_f \big) \in \Omega_{\delta}$ and 
\begin{equation}
\label{partie_re_positive1}
\Re \big( \lag \psi^j_f , \psi^{j,\eta}_{ref}(T^{\eta}) \rag \big) > 0, \quad \forall j \in \{1,2,3\}.
\end{equation}

\noindent
Let $u:= \Upsilon\big( \tilde{\psi}^1_f, \tilde{\psi}^2_f , \tilde{\psi}^3_f \big)$. Let $(\psi^1,\psi^2,\psi^3)$ be the solution of system~(\ref{syst_Neq})-(\ref{cond_ini}) with control $u$. We prove that
\begin{equation*}
\big( \psi^1(T^{\eta}), \psi^2(T^{\eta}), \psi^3(T^{\eta}) \big) = \big( \psi^1_f , \psi^2_f , \psi^3_f \big).
\end{equation*}
Up to a reduction of $\varepsilon_0$, we can assume that
\begin{equation}
\label{partie_re_positive2}
\Re \big( \lag \psi^j(T^{\eta}) , \psi^{j,\eta}_{ref}(T^{\eta}) \big) > 0, \quad \forall j \in \{1,2,3\}.
\end{equation}
By definition of $\Upsilon$ and $\tilde{\PP}_1$ it comes that
\begin{equation*}
\psi^1(T^{\eta}) - \Re (\lag \psi^1(T^{\eta}) , \psi^{1,\eta}_{ref}(T^{\eta}) \rag ) \psi^{1,\eta}_{ref}(T^{\eta})
=
\psi^1_f - \Re (\lag \psi^1_f , \psi^{1,\eta}_{ref}(T^{\eta}) \rag ) \psi^{1,\eta}_{ref}(T^{\eta}).
\end{equation*}
Thanks to (\ref{partie_re_positive1})-(\ref{partie_re_positive2}) and the fact that $|| \psi^1(T^{\eta}) ||_{L^2} = || \psi^1_f ||_{L^2}$, we get
\begin{equation}
\label{controle_NL_1}
\psi^1(T^{\eta}) = \psi^1_f.
\end{equation}

\noindent
The equality $\tilde{\PP}_2(\psi^2(T^{\eta}))=\tilde{\psi}^2_f$ gives
\begin{equation}
\label{controle_NL_2_0}
\begin{aligned}
&\psi^2(T^{\eta}) - \lag \psi^2(T^{\eta}),\psi^{1,\eta}_{ref}(T^{\eta}) \rag \psi^{1,\eta}_{ref}(T^{\eta}) - \Re (\lag \psi^2(T^{\eta}) , \psi^{2,\eta}_{ref}(T^{\eta}) \rag) \psi^{2,\eta}_{ref}(T^{\eta}) 
\\
&= \psi^2_f - \lag \psi^2_f,\psi^{1,\eta}_{ref}(T^{\eta}) \rag \psi^{1,\eta}_{ref}(T^{\eta}) - \Re (\lag \psi^2_f , \psi^{2,\eta}_{ref}(T^{\eta}) \rag) \psi^{2,\eta}_{ref}(T^{\eta}).
\end{aligned}
\end{equation}
Taking the scalar product of (\ref{controle_NL_2_0}) with $\psi^1_f$, using (\ref{controle_NL_1}) and the constraints $ \lag \psi^2(T^{\eta}),\psi^1(T^{\eta}) \rag = \lag \psi^2_f,\psi^1_f \rag = 0$, it comes that
\begin{equation}
\label{contraintes1}
\begin{aligned}
\lag \psi^2(T^{\eta}),\psi^{1,\eta}_{ref}(T^{\eta}) \rag \lag \psi^{1,\eta}_{ref}(T^{\eta}),\psi^1_f \rag + \Re \big( \lag \psi^2(T^{\eta}),\psi^{2,\eta}_{ref}(T^{\eta}) \rag \big) \lag \psi^{2,\eta}_{ref}(T^{\eta}),\psi^1_f \rag
\\
= \lag \psi^2_f,\psi^{1,\eta}_{ref}(T^{\eta}) \rag \lag \psi^{1,\eta}_{ref}(T^{\eta}),\psi^1_f \rag + \Re \big( \lag \psi^2_f,\psi^{2,\eta}_{ref}(T^{\eta}) \rag \big) \lag \psi^{2,\eta}_{ref}(T^{\eta}),\psi^1_f \rag.
\end{aligned}
\end{equation}
As $|| \psi^2(T^{\eta}) ||_{L^2} = || \psi^2_f ||_{L^2}$, we also get
\begin{equation}
\label{contraintes2}
\begin{aligned}
| \lag \psi^2(T^{\eta}),\psi^{1,\eta}_{ref}(T^{\eta}) \rag|^2 + \Re \big( \lag \psi^2(T^{\eta}),\psi^{2,\eta}_{ref}(T^{\eta}) \rag \big)^2 
\\
= | \lag \psi^2_f,\psi^{1,\eta}_{ref}(T^{\eta}) \rag|^2 + \Re \big( \lag \psi^2_f,\psi^{2,\eta}_{ref}(T^{\eta}) \rag \big)^2 .
\end{aligned}
\end{equation}
Straightforward computations prove that, up to an a priori reduction of $\varepsilon_0$, equalities (\ref{contraintes1}) and (\ref{contraintes2}) imply
\begin{equation}
\Re \big( \lag \psi^2(T^{\eta}),\psi^{2,\eta}_{ref}(T^{\eta}) \rag \big) = \Re \big( \lag \psi^2_f,\psi^{2,\eta}_{ref}(T^{\eta}) \rag \big)
\end{equation}
Then (\ref{contraintes1}), imply
$\lag \psi^2(T^{\eta}),\psi^{1,\eta}_{ref}(T^{\eta}) \rag = \lag \psi^2_f,\psi^{1,\eta}_{ref}(T^{\eta}) \rag$.
Finally, using these two last equalities in (\ref{controle_NL_2_0}), we obtain
\begin{equation}
\label{controle_NL_2}
\psi^2(T^{\eta}) = \psi^2_f.
\end{equation}

\noindent
Using $\tilde{\PP}_3(\psi^3(T^{\eta})) = \tilde{\psi}^3_f$ and the exact same strategy we also get
\begin{equation}
\label{controle_NL_3}
\psi^3(T^{\eta}) = \psi^3_f.
\end{equation}
Thus equalities (\ref{controle_NL_1}), (\ref{controle_NL_2}) and (\ref{controle_NL_3}) end the proof of Theorem~\ref{theo_controle_3eq} with $T^*:=T^\eta$ and
\begin{equation*}
\Gamma  : \big(\psi^1_f , \psi^2_f , \psi^3_f \big) \mapsto \Upsilon \big( \tilde{\PP}_1(\psi^1_f) , \tilde{\PP}_2(\psi^2_f) , \tilde{\PP}_3(\psi^3_f) \big).
\end{equation*}
\hfill $\blacksquare$
\\

\begin{rmq}
\label{rmq_th_cond_ini}
As mentioned in Remark~\ref{rmq_gestion_cond_ini}, a slight change in the proof allows to prove Theorem~\ref{theo_controle_3eq} for initial conditions $(\psi^1_0, \psi^2_0, \psi^3_0)$ close enough to $(\varphi_1,\varphi_2,\varphi_3)$ satisfying
\begin{equation}
\label{cond_contrainte}
\lag \psi^j_0 , \psi^k_0 \rag = \lag \psi^j_f , \psi^k_f \rag, \quad \forall j,k \in \{1,2,3\}.
\end{equation}
To this aim, the inverse mapping theorem is applied at the point $(u^{\eta}_{ref},\varphi_1,\varphi_2,\varphi_3)$ to the map
\begin{equation*}
\Lambda : L^2((0,T^{\eta}),\R)\times ( \mathcal{S} \cap H^3_{(0)}(0,1))^3 \to ( \mathcal{S} \cap H^3_{(0)}(0,1))^3 \times X^f_{T^{\eta}}
\end{equation*}
defined by
\begin{equation*}
\Lambda   \big( u, \psi^1_0 , \psi^2_0 , \psi^3_0 \big) =
\big( (\psi^j_0)_{j=1,2,3} , \tilde{\PP}_j(\psi^j(T^{\eta}))_{j=1,2,3} \big).
\end{equation*}
The compatibility condition (\ref{cond_contrainte}) will then lead to (\ref{contraintes1}), the conclusion being unchanged.
\end{rmq}

\section{Controllability results for two equations}
\label{sect_2eq}

\noindent
Theorem~\ref{theo_controle_3eq} leads to local exact controllability up to a global phase and a global delay in the case $N=2$. Actually the strategy we developed can be improved in this case to obtain less restrictive results, namely Theorems \ref{th_controle_2eq_1} and \ref{th_controle_2eq_2}.
Here, we only detail the construction of the reference trajectory, the application of the return method being very similar to Section \ref{sect_methode_retour}. Subsection \ref{subsect_controle_2eq_1} will imply Theorem~\ref{th_controle_2eq_1} and Subsection \ref{subsect_controle_2eq_2} will imply Theorem~\ref{th_controle_2eq_2}.

\noindent
In all this section, we consider $N=2$. Let $T_1>0$ and $\varepsilon \in (0,T_1)$. As in Theorem~\ref{th_construction_uref}, the reference control is designed in two steps.

\noindent
Let $u \equiv 0$ on $[0, \frac{\varepsilon}{2})$. Proposition~\ref{conditions_minimalite} is replaced by the following proposition.
\begin{prop}
\label{conditions_minimalite_2eq}
There exists $\eta^* >0$ and a $C^1$ map
\begin{equation*}
\hat{\Gamma} : (0,\eta^*) \rightarrow L^2 \left( \left(\frac{\varepsilon}{2},\varepsilon \right),\R \right),
\end{equation*}
satisfying $\hat{\Gamma}(0)= 0$ such that for any $\eta \in (0,\eta^*)$, the solution $(\psi^{1,\eta}_{ref},\psi^{2,\eta}_{ref})$ of system~(\ref{syst_Neq}) with control $u:= \hat{\Gamma}(\eta)$ and initial conditions $\psi^{j,\eta}_{ref}(\frac{\varepsilon}{2}) = \Phi_j(\frac{\varepsilon}{2})$ for $j=1,2$ satisfies
\begin{align*}
\lag \mu \psi^{1,\eta}_{ref} (\varepsilon) , \psi^{1,\eta}_{ref} (\varepsilon) \rag &= \lag \mu \varphi_1 , \varphi_1 \rag + \eta,
\\
\lag \mu \psi^{2,\eta}_{ref} (\varepsilon) , \psi^{2,\eta}_{ref} (\varepsilon) \rag &= \lag \mu \varphi_2 , \varphi_2 \rag.
\end{align*}
\end{prop}
As previously, this proposition will ensure controllability of the linearized system around the reference trajectory.
The proof is a simple adaptation of Proposition~\ref{conditions_minimalite} and is not detailed.

\noindent
We now turn to two different constructions of reference trajectories on $(\varepsilon,T_1)$, to replace Proposition~\ref{traj_ref_inversion_locale}.

\subsection{Controllability up to a global phase in arbitrary time : Theorem \ref{th_controle_2eq_1}}
\label{subsect_controle_2eq_1}

Let $T>0$ be arbitrary. Up to a reduction of $\varepsilon$, we assume that $T=T_1$. We prove that there exists a global phase $\theta^{\eta} >0$ and a control $u_{ref}^{\eta}$ on $(\varepsilon,T)$ such that the associated trajectory $(\psi^{1,\eta}_{ref},\psi^{2,\eta}_{ref})$ of (\ref{syst_Neq})-(\ref{cond_ini}) satisfies Proposition~\ref{conditions_minimalite_2eq}, 
\begin{equation}
\label{2eq_temps_final_1}
(\psi^{1,\eta}_{ref}(T) , \psi^{2,\eta}_{ref}(T)) = e^{i \theta^{\eta}} ( \Phi_1(T) , \Phi_2(T)  ),
\end{equation}
and $||u_{ref}^{\eta}||_{L^2(0,T)} \leq C \eta$.

\noindent
Proposition~\ref{traj_ref_inversion_locale} is replaced by the following proposition which proof is a simple adaptation of the one of Proposition~\ref{traj_ref_inversion_locale} and is not detailed.
\begin{prop}
\label{traj_ref_inversion_locale_2eq_1}
There exists $\delta > 0$ and a $C^1$-map 
\begin{equation*}
\tilde{\Gamma} : \tilde{\OO}_{\delta} \to L^2((\varepsilon,T) ,\R)
\end{equation*}
with 
\begin{equation*}
\tilde{\OO}_{\delta} := \left\{ \big( \psi^1_0,\psi^2_0 \big) \in (\mathcal{S} \cap H^3_{(0)}(0,1))^2 \: ; \: 
\sum_{j=1}^2 ||\psi^j_0 - \Phi_j(\varepsilon)||_{H^3_{(0)}} < \delta \right\},
\end{equation*}
such that $\tilde{\Gamma} \big( \Phi_1(\varepsilon), \Phi_2(\varepsilon) \big)=0$ and, if $(\psi^1_0,\psi^2_0) \in \tilde{\OO}_{\delta}$, the solution $\big( \psi^1, \psi^2 \big)$ of system~(\ref{syst_Neq}) with initial conditions $\psi^j(\varepsilon,\cdot)=\psi^j_0$, for $j =1,2$, and control $u:= \tilde{\Gamma} \big( \psi^1_0,\psi^2_0 \big)$ satisfies
\begin{align}
\label{uref_condition_2eq_1}
&\PP_1 \big( \psi^1(T) \big) = \PP_2 \big( \psi^2(T) \big) =  0,
\\
\label{uref_condition_2eq_2}
&\Im \big( \lag \psi^1(T) ,\Phi_1(T) \rag \overline{\lag \psi^2(T) , \Phi_2(T) \rag}  \big) \Big) = 0.
\end{align}
\end{prop}

\noindent
There exists $\overline{\eta} >0$ such that for $\eta \in (0, \overline{\eta})$, the control
\begin{equation}
\label{def_uref_2eq}
u_{ref}^{\eta} (t) :=
\left\{
\begin{aligned}
& 0 \quad & \text{for } t \in (0,\frac{\varepsilon}{2}),
\\
\hat{\Gamma} &(\eta) \quad & \text{for } t \in (\frac{\varepsilon}{2},\varepsilon),
\\
\tilde{\Gamma}(\psi^{1,\eta}_{ref}(\varepsilon), &\psi^{2,\eta}_{ref}(\varepsilon) ) \quad & \text{for } t \in (\varepsilon,T),
\end{aligned}
\right.
\end{equation}
is well defined and satisfies $||u_{ref}^{\eta}||_{L^2(0,T)} \leq C \eta$, where $\hat{\Gamma}$ and $\tilde{\Gamma}$ are defined respectively in Proposition \ref{conditions_minimalite_2eq} and \ref{traj_ref_inversion_locale_2eq_1}.
Proposition~\ref{traj_ref_inversion_locale_2eq_1} implies that
\begin{align*}
&\psi^{1,\eta}_{ref}(T) = \lag  \psi^{1,\eta}_{ref}(T) , \Phi_1(T) \rag \Phi_1(T),
\\
&\psi^{2,\eta}_{ref}(T) = \lag \psi^{2,\eta}_{ref}(T) , \Phi_1(T) \rag \Phi_1(T) + \lag \psi^{2,\eta}_{ref}(T) , \Phi_2(T) \rag \Phi_2(T),
\\
&\Im \big( \lag \psi^{1,\eta}_{ref}(T) ,\Phi_1(T) \rag \overline{\lag \psi^{2,\eta}_{ref}(T) , \Phi_2(T) \rag} \big) =0.
\end{align*}
Thus, using the invariant of the system, it comes that there exist $\theta^{\eta}_1$, $\theta^{\eta}_2 \in [0,2\pi)$ such that
\begin{equation*}
(\psi^{1,\eta}_{ref}(T) , \psi^{2,\eta}_{ref}(T)) = \big( e^{-i\theta^{\eta}_1} \Phi_1(T) , e^{-i\theta^{\eta}_2} \Phi_2(T) \big),
\end{equation*}
and
\begin{equation*}
\theta^{\eta}_1-\theta^{\eta}_2 \equiv 0 \, [2 \pi].
\end{equation*}
Finally, this implies that there exists $\theta^\eta \in \R$ such that
\begin{equation*}
(\psi^{1,\eta}_{ref}(T) , \psi^{2,\eta}_{ref}(T)) = e^{i \theta^{\eta}} ( \Phi_1(T) , \Phi_2(T)  ).
\end{equation*}
Then, application of the return method along this trajectory as in Section \ref{sect_methode_retour} implies Theorem~\ref{th_controle_2eq_1}.

\begin{rmq}
To investigate controllability properties up to a global phase, as proposed in \cite{MRT}, one can introduce a fictitious control $\omega$ in the following way
\begin{equation*}
\left\{
\begin{aligned}
& i \partial_t \psi^j = -\partial^2_{xx} \psi^j - u(t) \mu (x) \psi^j - \omega(t) \psi^j, \: &(t,x)& \in (0,T) \times (0,1), \; j \in \{1,2\},
\\
& \psi^j(t,0) = \psi^j(t,1) =0, \: &t& \in (0,T), \; j \in \{1,2\}.
\end{aligned}
\right.
\end{equation*}
Adapting the strategy of \cite[Theorem 1]{BeauchardLaurent}, one can prove local controllability of this system by linearization around the trajectory $(\Phi_1, \Phi_2, u\equiv0, \omega \equiv 0)$. This would lead to local controllability up to a global phase. However, in this case, one would obtain for each target $(\psi^1_f,\psi^2_f)$ close enough to $(\Phi_1,\Phi_2)$ a global phase $\theta=\theta(\psi^1_f,\psi^2_f)$ such that there exists a control driving the solution of (\ref{syst_Neq}) from (\ref{cond_ini}) to $e^{i \theta} (\psi^1_f,\psi^2_f)$. 
\end{rmq}

\subsection{Exact controllability up to a global delay : Theorem \ref{th_controle_2eq_2}}
\label{subsect_controle_2eq_2}

We prove that there exists $T^{\eta}>0$ and a control $u_{ref}^{\eta}$ on $(\varepsilon,T_1)$ such that if $u_{ref}^{\eta}$ is extended by $0$ on $(T_1,T^{\eta})$, the associated trajectory $(\psi^{1,\eta}_{ref},\psi^{2,\eta}_{ref})$ of (\ref{syst_Neq})-(\ref{cond_ini}) satisfies Proposition~\ref{conditions_minimalite_2eq}, 
\begin{equation}
\label{2eq_temps_final_2}
(\psi^{1,\eta}_{ref}(T^{\eta}) , \psi^{2,\eta}_{ref}(T^{\eta})) =  ( \varphi_1 , \varphi_2  ),
\end{equation}
and $||u_{ref}^{\eta}||_{L^2(0,T^{\eta})} \leq C \eta$.

\noindent
Proposition~\ref{traj_ref_inversion_locale} is replaced by the following proposition which proof is a simple adaptation of the one of Proposition~\ref{traj_ref_inversion_locale} and is not detailed.
\begin{prop}
\label{traj_ref_inversion_locale_2eq_2}
There exists $\delta > 0$ and a $C^1$-map 
\begin{equation*}
\tilde{\Gamma} : \tilde{\OO}_{\delta} \to L^2((\varepsilon,T_1) ,\R)
\end{equation*}
with 
\begin{equation*}
\tilde{\OO}_{\delta} := \left\{ \big( \psi^1_0,\psi^2_0 \big) \in (\mathcal{S} \cap H^3_{(0)}(0,1))^2 \: ; \: 
\sum_{j=1}^2 ||\psi^j_0 - \Phi_j(\varepsilon)||_{H^3_{(0)}} < \delta \right\},
\end{equation*}
such that $\tilde{\Gamma} \big( \Phi_1(\varepsilon), \Phi_2(\varepsilon) \big)=0$ and, if $(\psi^1_0,\psi^2_0) \in \tilde{\OO}_{\delta}$, the solution $\big( \psi^1, \psi^2 \big)$ of system~(\ref{syst_Neq}) with initial conditions $\psi^j(\varepsilon,\cdot)=\psi^j_0$, for $j =1,2$, and control $u:= \tilde{\Gamma} \big( \psi^1_0,\psi^2_0 \big)$ satisfies
\begin{align}
\label{uref_condition_2eq_3}
&\PP_1 \big( \psi^1(T_1) \big) = \PP_2 \big( \psi^2(T_1) \big) =  0,
\\
\label{uref_condition_2eq_4}
&\Im \Big( \lag \psi^1(T_1) ,\Phi_1(T_1) \rag^4 \overline{\lag \psi^2(T_1) , \Phi_2(T_1) \rag}  \big) \Big) = 0.
\end{align}
\end{prop}
\noindent
There exists $\overline{\eta} >0$ such that for $\eta \in (0, \overline{\eta})$, the control
\begin{equation}
\label{def_uref_2eq2}
u_{ref}^{\eta} (t) :=
\left\{
\begin{aligned}
& 0 \quad & \text{for } t \in (0,\frac{\varepsilon}{2}),
\\
\hat{\Gamma} &(\eta) \quad & \text{for } t \in (\frac{\varepsilon}{2},\varepsilon),
\\
\tilde{\Gamma}(\psi^{1,\eta}_{ref}(\varepsilon), &\psi^{2,\eta}_{ref}(\varepsilon) ) \quad & \text{for } t \in (\varepsilon,T_1),
\end{aligned}
\right.
\end{equation}
is well defined and satisfies $||u_{ref}^{\eta}||_{L^2(0,T_1)} \leq C \eta$, where $\hat{\Gamma}$ and $\tilde{\Gamma}$ are defined respectively in Proposition \ref{conditions_minimalite_2eq} and \ref{traj_ref_inversion_locale_2eq_2}. 
Proposition~\ref{traj_ref_inversion_locale_2eq_2} implies the existence of $\theta^{\eta}_1$, $\theta^{\eta}_2 \in [0, 2 \pi)$ such that 
\begin{align*}
(\psi^{1,\eta}_{ref}(T_1) , \psi^{2,\eta}_{ref}(T_1)) &= \big( e^{-i\theta^{\eta}_1} \Phi_1(T_1) , e^{-i\theta^{\eta}_2} \Phi_2(T_1) \big),
\\
4 \theta^{\eta}_1 - \theta^{\eta}_2 \, &\equiv \, 0 \, [2\pi].
\end{align*}
Let $T^{\eta} > T_1$ be such that
\begin{equation*}
\theta^{\eta}_1 + \lambda_1 T^{\eta} \, \equiv \, 0 \, [2 \pi]
\end{equation*}
Thus,
\begin{equation*}
\theta^{\eta}_2 + \lambda_2 T^{\eta} \, \equiv \, 4 \big( \theta^{\eta}_1 + \lambda_1 T^{\eta} \big) \, \equiv \, 0 \, [2 \pi].
\end{equation*}
Finally, if we extend $u_{ref}^{\eta}$ by $0$ on $(T_1,T^{\eta})$, we have that $(\psi^{1,\eta}_{ref},\psi^{2,\eta}_{ref})$ is solution of (\ref{syst_Neq})-(\ref{cond_ini}) with control $u_{ref}^{\eta}$ and satisfies
\begin{equation*}
\psi^{j,\eta}_{ref}(T^{\eta}) = e^{-i( \theta^{\eta}_j + \lambda_j T^{\eta} )} \varphi_j = \varphi_j.
\end{equation*}
Then, application of the return method along this trajectory as in Section \ref{sect_methode_retour} implies Theorem~\ref{th_controle_2eq_2}.

\section{Non controllability results in small time}
\label{sect_non_controlabilite}

The goal of this section is the proof of Theorems \ref{th_NC_2eq_1} and \ref{th_NC_3eq}.

\subsection{Heuristic of non controllability}

\noindent
We adapt the strategy developed in \cite{BM_Tmin} by Beauchard and the author in the case $N=1$. Using power series expansion, we consider 
\begin{equation} \label{power_series_expansion}
\begin{aligned}
u &= 0 + \varepsilon v,
\\
\psi^j &= \Phi_j + \varepsilon \Psi^j + \varepsilon^2 \xi^j + o(\varepsilon^2), \quad \forall j \in \{1,\dots,N\}.
\end{aligned}
\end{equation}
Here and in the following, we use the classical Landau notations. We say that $f = \underset{x \to a}{O}(g)$ if there exist $C>0$ and a neighbourhood $\VV(a)$ of $a$ such that $||f(x)||\leq C||g(x)||$ for $x \in \VV(a)$. We say that $f = \underset{x \to a}{o}(g)$ if for any $\delta>0$ there exists a neighbourhood $\VV(a)$ of $a$ such that $||f(x)||\leq \delta||g(x)||$ for $x \in \VV(a)$.

\noindent
Considering (\ref{power_series_expansion}), we define the following systems for $j \in \{1,\dots,N\}$,
\begin{equation}
\label{2eq_ordre1}
\left\{
\begin{aligned}
&i \partial_t \Psi^j = - \partial^2_{xx} \Psi^j - v(t) \mu(x) \Phi_j, \: &(t,x)& \in (0,T) \times (0,1),
\\
&\Psi^j(t,0) = \Psi^j(t,1) = 0, \: &t& \in (0,T),
\\
&\Psi^j(0,x) = 0, \: &x& \in (0,1),
\end{aligned}
\right.
\end{equation}
and
\begin{equation}
\label{2eq_ordre2}
\left\{
\begin{aligned}
&i \partial_t \xi^j = - \partial^2_{xx} \xi^j - v(t) \mu(x) \Psi_j - w(t) \mu(x) \Phi_j, \, &(t,x)& \in (0,T) \times (0,1),
\\
&\xi^j(t,0) = \xi^j(t,1) = 0, \: &t& \in(0,T),
\\
&\xi^j(0,x) = 0, \: &x& \in (0,1).
\end{aligned}
\right.
\end{equation}
We focus in this heuristic on the case $N=2$. Let us try to reach
\begin{equation}
\label{cible_2eq}
\big( \psi^1(T), \psi^2(T) \big) = \Big( \Phi_1(T) , (\sqrt{1-\delta^2} + i \alpha \delta) \Phi_2(T) \Big),
\end{equation}
with $\delta>0$ and $\alpha$ defined in Theorem~\ref{th_NC_2eq_1} from $(\psi^1(0),\psi^2(0)) = (\varphi_1, \varphi_2)$. 
Condition (\ref{cible_2eq}) imposes $\Psi^1(T) = 0$ i.e.
\begin{equation*}
v \in V_T := \left\{ v \in L^2((0,T),\R) \, ; \, \int_0^T v(t) e^{i (\lambda_k - \lambda_1)t} \md t =0 , \: \forall k \in \N^* \right\}.
\end{equation*}

\noindent
Let us define the following quadratic forms, for $j \in \{1,2\}$, associated to the second order
\begin{align*}
Q_{T,j}(v) :&= \Im \big( \lag \xi^j(T) , \Phi_j(T) \rag \big)
\\
&= \int_0^T v(t) \int_0^t v(\tau) \left( \sum_{k=1}^{+\infty} \lag \mu \varphi_j, \varphi_k \rag^2
\sin ((\lambda_k - \lambda_j)(t-\tau)) \right) \md \tau \md t,
\end{align*}
and
\begin{equation}
\label{def_fq}
Q_T(v) := \lag \mu \varphi_1,\varphi_1 \rag Q_{T,2}(v) - \lag \mu \varphi_2,\varphi_2 \rag Q_{T,1}(v).
\end{equation}
The following proposition states that in time small enough, the quadratic form $Q_T$ has a sign on $V_T$.

\begin{prop}
\label{prop_heuristique_NC}
Assume that $\mu$ satisfies Hypothesis~\ref{hypo_mu2}. Then, there exists $T_*>0$ such that for any $T \in (0, T_*)$, for any $v \in V_T \backslash \{0\}$, 
\begin{equation*}
\AA Q_T(v) < 0,
\end{equation*}
where $\AA \in \R^*$ is defined in Hypothesis~\ref{hypo_mu2}.
\end{prop}

\noindent
\textbf{Proof of Proposition~\ref{prop_heuristique_NC} :}
Let $v \in V_T$ and $s : t \in (0,T) \mapsto \int_0^t v(\tau) \md \tau$. Performing integrations by part, we define a new quadratic form
\begin{equation}
\label{def_fq_js}
\QQ_{T,j} (s) := -\lag(\mu')^2 \varphi_j ,\varphi_j \rag \int_0^T s(t)^2 \md t + \int_0^T s(t) \int_0^t s(\tau) h_j(t-\tau) \md \tau \md t = Q_{T,j}(v),
\end{equation}
where $h_j : t \mapsto \sum_{k=1}^{+\infty} (\lambda_k - \lambda_j)^2 \lag \mu \varphi_j, \varphi_k \rag^2 \sin ((\lambda_k - \lambda_j) t)$. As $\mu \in H^3((0,1),\R)$, it comes that $h_j\in C^0(\R,\R)$.
Thus, if we define
\begin{equation}
\label{def_fq_s}
\QQ_T(s) := \lag \mu \varphi_1 , \varphi_1 \rag \QQ_{T,2}(s) - \lag \mu \varphi_2 ,\varphi_2 \rag \QQ_{T,1}(s),
\end{equation}
we get that
\begin{equation*}
Q_T(v) = \QQ_T(s) = -\AA ||s||_{L^2}^2 + \int_0^T s(t) \int_0^t s(\tau) h(t-\tau) \md \tau \md t,
\end{equation*}
with
\begin{equation*}
h := \lag \mu \varphi_1 , \varphi_1 \rag h_2 - \lag \mu \varphi_2 ,\varphi_2 \rag h_1 \in C^0(\R,\R).
\end{equation*}
We can assume, without loss of generality, that $\AA>0$. Thus, there exists $C=C(\mu)>0$ such that
\begin{equation}
\label{coercivite}
Q_T(v) \leq \big( -\AA + CT \big) ||s||_{L^2}^2.
\end{equation}
We conclude the proof by choosing $T_* < \frac{\AA}{C}$. 

\hfill $\blacksquare$

\begin{rmq}
This Proposition indicates that, in small time, there are targets that cannot be reached. However, using the theory of Legendre form (see e.g. \cite{Hestenes, BonnansShapiro}), we can prove that $Q_T$ lacks coercivity in $L^2((0,T),\R)$. This is why we work directly with the quadratic form $\QQ_T$ adapted to the auxiliary system defined in Subsection \ref{subsect_syst_aux} where the control is $s$ and not $v$.
\end{rmq}

\begin{rmq}
This strategy is only valid for small time and we do not know if this quadratic form changes sign in time large enough on $V_T$. Following the strategy of \cite{BM_Tmin}, this would imply local exact controllability in large time but it is an open question.
\end{rmq}

\subsection{Auxiliary system}
\label{subsect_syst_aux}

\noindent
For $j \in \{ 1, \dots ,N \}$ , we consider the function $\tilde{\psi}^j$ defined by
\begin{equation}
\label{chgt_var}
\psi^j(t,x) = \tilde{\psi}^j(t,x) e^{i s(t) \mu(x)}
\, \text{ with } 
s(t) := \int_0^t u(\tau) \md \tau.
\end{equation}
It is a weak solution of
\begin{equation}
\label{syst_aux}
\left\{
\begin{aligned}
&i \partial_t \tilde{\psi}^j = -\partial^2_{xx} \tilde{\psi}^j -i s(t) \big( 2 \mu'(x) \partial_x \tilde{\psi}^j + \mu''(x) \tilde{\psi}^j \big) + s(t)^2 \mu'(x)^2 \tilde{\psi}^j, 
\\
&\tilde{\psi}^j(t,0) = \tilde{\psi}^j(t,1) = 0,
\\
&\tilde{\psi}^j (0,  \cdot) = \varphi_j.
\end{aligned}
\right.
\end{equation}
Using Proposition~\ref{bien_pose} on (\ref{syst_Neq}) and (\ref{chgt_var}), it follows that the following well posedness result holds. In the following, the time derivative of $s$ will be denoted by $\dot{s}$.

\begin{prop}
\label{bien_pose_aux}
Let $\mu \in H^3((0,1),\R)$, $T>0$, $s \in H^1((0,T),\R)$ with $s(0)=0$.
There exists a unique weak solution $(\tilde{\psi}^1, \dots, \tilde{\psi}^N) \in C^0([0,T],H^3 \cap H^1_0)^N$ of system~(\ref{syst_aux}).
Moreover, for every $R>0$, there exists $C=C(T,\mu,R)>0$ such that,
if $\|\dot{s}\|_{L^2(0,T)} < R$, then this weak solution satisfies for any $j \in \{1, \dots, N\}$,
\begin{equation*}
\| \tilde{\psi}^j \|_{L^\infty((0,T),H^3 \cap H^1_0)} \leqslant C.
\end{equation*}
\end{prop}

\subsection{Non exact controllability in arbitrary time with $N=2$.}

In this subsection, we consider system~(\ref{syst_Neq}) with $N=2$ and prove Theorem~\ref{th_NC_2eq_1}. 
This result is a corollary of the following theorem for the auxiliary system.
\begin{theo}
\label{th_NC_aux_2eq_1}
Let $\mu \in H^3((0,1),\R)$ be such that Hypothesis~\ref{hypo_mu2} hold. Let $T_*>0$ be as in Proposition~\ref{prop_heuristique_NC} and $\alpha \in \{ -1, 1 \}$ as in Theorem~\ref{th_NC_2eq_1}. For any $T<T_*$, there exists $\varepsilon>0$ such that for every $s \in H^1((0,T),\R)$ with $s(0)=0$ and $||\dot{s}||_{L^2} < \varepsilon$, the solution of system~(\ref{syst_aux}) satisfies
\begin{equation*}
\big( \tilde{\psi}^1(T), \tilde{\psi}^2(T) \big) \neq \left( \Phi_1(T) e^{i \theta \mu} ,
\left( \sqrt{1-\delta^2} + i \alpha \delta \right) \Phi_2(T) e^{i \theta \mu} \right), \, \forall \delta >0 ,
\forall \theta \in \R.
\end{equation*}
\end{theo}
Before getting into the proof of Theorem~\ref{th_NC_aux_2eq_1}, we prove that it implies Theorem~\ref{th_NC_2eq_1}.

\noindent
\textbf{Proof of Theorem~\ref{th_NC_2eq_1} : }
Let $T<T_*$ and $\varepsilon>0$ defined by Theorem~\ref{th_NC_aux_2eq_1}. Let $u \in L^2((0,T),\R)$ be such that 
$|| u ||_{L^2(0,T)} < \varepsilon$.
Assume by contradiction that
\begin{equation*}
\big( \psi^1(T), \psi^2(T) \big) = \left( \Phi_1(T) ,
\left( \sqrt{1-\delta^2} + i \alpha \delta \right) \Phi_2(T) \right), 
\end{equation*}
for some $\delta >0$. Let $s$ and $\tilde{\psi}^j$ be defined by (\ref{chgt_var}). Then $s(0)=0$, $||\dot{s}||_{L^2}< \varepsilon$ and $\tilde{\psi}^j$ is solution of (\ref{syst_aux}) and satisfies
\begin{equation*}
\big( \tilde{\psi}^1(T), \tilde{\psi}^2(T) \big) = \left( \Phi_1(T) e^{-i s(T) \mu} ,
\left( \sqrt{1-\delta^2} + i \alpha \delta \right) \Phi_2(T) e^{-i s(T) \mu} \right).
\end{equation*}
Thanks to Theorem~\ref{th_NC_aux_2eq_1}, this is impossible.

\hfill $\blacksquare$

\noindent
\textbf{Proof of Theorem~\ref{th_NC_aux_2eq_1} :}
Without loss of generality, we assume that $\AA>0$.
\\

\textit{First step :} we prove that $-\QQ_T$ is coercive for $T<T_*$.

\noindent
Using the same estimates as in (\ref{coercivite}) and the fact that $T_* < \frac{\AA}{C}$, we get that there exists $C_*>0$ such that for $T<T_*$
\begin{equation}
\label{coercivite2}
\QQ_T(s) \leq -C_* ||s||_{L^2}^2, \quad \forall s \in L^2((0,T),\R).
\end{equation}
\\

\textit{Second step :} approximation of first and second order.

\noindent
Using the first and second order approximation of (\ref{syst_aux}), the following lemma holds.
\begin{lemme}
\label{lemme_approximation_aux}
Let $T>0$ and $\mu \in H^3((0,1),\R)$. For all $j \in \{ 1, \dots , N\}$
\begin{align*}
\left| \Im ( \lag \tilde{\psi}^j(T) ,\Phi_j(T) \rag ) - \QQ_{T,j}(s) \right| = o(||s||_{L^2}^2) 
\text{ when } ||\dot{s} ||_{L^2} \to 0,
\\
\left| \Im ( \lag \tilde{\psi}^j(T) ,\Phi_j(T) \rag ) \right| = o(||s||_{L^2}) 
\text{ when } ||\dot{s} ||_{L^2} \to 0.
\end{align*}
\end{lemme}

\noindent
\textbf{Proof of Lemma~\ref{lemme_approximation_aux}.} Let $j \in \{1, \dots, N\}$. As proved in \cite[Proposition 3]{BM_Tmin}, if we define the first and second order approximations, $\tilde{\Psi}^j$ and $\tilde{\xi}^j$, by
\begin{equation} \label{ordre1_aux}
\Psi^j(t,x) = \tilde{\Psi}^j(t,x) + i s(t) \mu(x) \Phi_j(t,x),
\end{equation}
and
\begin{equation} \label{ordre2_aux}
\xi^j(t,x) = \tilde{\xi}^j(t,x) + i s(t) \mu(x) \tilde{\Psi}^j(t,x) - s(t)^2 \mu'(x)^2 \Phi_j(t,x),
\end{equation}
it comes that, when $||\dot{s}||_{L^2} \to 0$
\begin{equation} \label{approx_ordre1_aux}
||\tilde{\psi}^j - \Phi_j - \tilde{\Psi}^j ||_{L^\infty((0,T),H^1_0)} = o(||s||_{L^2})
\end{equation}
and
\begin{equation} \label{approx_ordre2_aux}
||\tilde{\psi}^j - \Phi_j - \tilde{\Psi}^j - \tilde{\xi}^j ||_{L^\infty((0,T),L^2)} = o(||s||_{L^2}^2).
\end{equation}
Straightforward computations using (\ref{ordre1_aux}) imply $\Im(\lag \tilde{\Psi}^j(T), \Phi_j(T) \rag) = 0$. Thus, from (\ref{approx_ordre1_aux}) we deduce
\begin{equation*}
\left| \Im ( \lag \tilde{\psi}^j(T) ,\Phi_j(T) \rag ) \right| =
\left| \Im ( \lag (\tilde{\psi}^j - \Phi_j - \tilde{\Psi}^j)(T) ,\Phi_j(T) \rag ) \right| =
\underset{||\dot{s} ||_{L^2} \to 0}{o}(||s||_{L^2}).
\end{equation*}
Straightforward computations using (\ref{ordre2_aux}) imply $\Im(\lag \tilde{\xi}^j(T), \Phi_j(T) \rag) = \QQ_{T,j}(s)$. Thus, from (\ref{approx_ordre2_aux}) we deduce
\begin{align*}
\left| \Im ( \lag \tilde{\psi}^j(T) ,\Phi_j(T) \rag ) - \QQ_{T,j}(s) \right| 
&= \left| \Im ( \lag (\tilde{\psi}^j - \Phi_j - \tilde{\Psi}^j - \tilde{\xi}^j)(T) ,\Phi_j(T) \rag ) \right|
\\
&=\underset{||\dot{s} ||_{L^2} \to 0}{o}(||s||_{L^2}^2).
\end{align*}
This ends the proof of Lemma~\ref{lemme_approximation_aux}.
\hfill $\blacksquare$
\\

\noindent
\textit{Third step :} conclusion.

\noindent
Let $T<T_*$. Assume by contradiction, that $\forall \varepsilon>0$, $\exists s_{\varepsilon} \in H^1((0,T),\R)$ with $s_\varepsilon(0)=0$ and $||\dot{s}_{\varepsilon}||_{L^2} < \varepsilon$ such that the associated solution of (\ref{syst_aux}) satisfies
\begin{equation*}
\big( \tilde{\psi}^1_{\varepsilon}(T), \tilde{\psi}^2_{\varepsilon}(T) \big) = \left( \Phi_1(T) e^{i \theta_{\varepsilon} \mu} ,
\left( \sqrt{1-\delta_{\varepsilon}^2} + i \alpha \delta_{\varepsilon} \right) \Phi_2(T) e^{i \theta_{\varepsilon} \mu} \right),
\end{equation*}
with $\delta_{\varepsilon} >0$ and $\theta_{\varepsilon} \in \R$. Notice that 
\begin{equation*}
\delta_{\varepsilon} \underset{\varepsilon \to 0}{\rightarrow} 0, \quad \theta_{\varepsilon} \underset{\varepsilon \to 0}{\rightarrow} 0.
\end{equation*}
Explicit computations lead to
\begin{equation*}
\Im \big(\lag \tilde{\psi}^1_{\varepsilon}(T) , \Phi_1(T) \rag \big) 
= \lag \mu \varphi_1, \varphi_1 \rag \theta_{\varepsilon} + \underset{\varepsilon \to 0}{O}(\theta_{\varepsilon}^3),
\end{equation*}
and
\begin{equation*}
\Im \big(\lag \tilde{\psi}^2_{\varepsilon}(T) , \Phi_2(T) \rag \big)
= \alpha \delta_{\varepsilon} + \sqrt{1-\delta_{\varepsilon}^2} \lag \mu \varphi_2 , \varphi_2 \rag \theta_{\varepsilon} 
+ \underset{\varepsilon \to 0}{O}(\theta_{\varepsilon}^2).
\end{equation*}
Thus, it comes that 
\begin{align*}
&\lag \mu \varphi_1,\varphi_1 \rag \Im ( \lag \tilde{\psi}^2_{\varepsilon}(T) , \Phi_2(T) \rag)  - \lag \mu \varphi_2, \varphi_2 \rag \Im (\lag \tilde{\psi}^1_{\varepsilon}(T) , \Phi_1(T) \rag)
\\
&=  \alpha \lag \mu \varphi_1,\varphi_1 \rag \delta_{\varepsilon}
- \lag \mu \varphi_1 ,\varphi_1 \rag \lag \mu \varphi_2 , \varphi_2 \rag \frac{\delta_{\varepsilon}^2}{\sqrt{1-\delta_{\varepsilon}^2}+1} \theta_{\varepsilon}
+ \underset{\varepsilon \to 0}{O}(\theta_{\varepsilon}^2).
\end{align*}
Using Lemma~\ref{lemme_approximation_aux} to estimate $\Im \big(\lag \tilde{\psi}^1_{\varepsilon}(T) , \Phi_1(T) \rag \big) $ and $\Im \big(\lag \tilde{\psi}^2_{\varepsilon}(T) , \Phi_2(T) \rag \big) $ it comes that 
\begin{equation*}
\theta_{\varepsilon} = \underset{\varepsilon \to 0}{o}(||s_{\varepsilon}||_{L^2}),
\quad
\delta_{\varepsilon} = \underset{\varepsilon \to 0}{o}(||s_{\varepsilon}||_{L^2}).
\end{equation*}
Thus,
\begin{align*}
\lag & \mu \varphi_1,\varphi_1 \rag \Im ( \lag \tilde{\psi}^2_{\varepsilon}(T) , \Phi_2(T) \rag)  - \lag \mu \varphi_2, \varphi_2 \rag \Im (\lag \tilde{\psi}^1_{\varepsilon}(T) , \Phi_1(T) \rag)
\\
&= \lag \mu \varphi_1,\varphi_1 \rag \alpha \delta_{\varepsilon} + \underset{\varepsilon \to 0}{o}(||s_{\varepsilon}||_{L^2}^2).
\end{align*}
Finally, combining this with Lemma~\ref{lemme_approximation_aux} and (\ref{coercivite2}), we obtain
\begin{align*}
0 &<  \alpha \lag \mu \varphi_1,\varphi_1 \rag \delta_{\varepsilon}
\\
&= \lag \mu \varphi_1,\varphi_1 \rag \Im ( \lag \tilde{\psi}^2_{\varepsilon}(T) , \Phi_2(T) \rag)  - \lag \mu \varphi_2, \varphi_2 \rag \Im (\lag \tilde{\psi}^1_{\varepsilon}(T) , \Phi_1(T) \rag) + \underset{\varepsilon \to 0}{o}(||s_{\varepsilon}||_{L^2}^2)
\\
&= \QQ_T(s_{\varepsilon}) + \underset{\varepsilon \to 0}{o}(||s_{\varepsilon}||_{L^2}^2)
\\
& \leq -C_* ||s_{\varepsilon}||_{L^2}^2 + \underset{\varepsilon \to 0}{o}(||s_{\varepsilon}||_{L^2}^2).
\end{align*}
This is impossible for $\varepsilon$ sufficiently small. This ends the proof of Theorem~\ref{th_NC_aux_2eq_1}.

\hfill $\blacksquare$

\subsection{Non exact controllability up to a global phase in arbitrary time with $N=3$.}

\noindent
In this subsection, we consider system~(\ref{syst_Neq}) with $N=3$ and prove Theorem~\ref{th_NC_3eq}. 
As previously, this result is a corollary of the following theorem for the auxiliary system.
\begin{theo}
\label{th_NC_aux_3eq}
Let $\mu \in H^3((0,1),\R)$ be such that Hypothesis~\ref{hypo_mu3} hold. Let $\beta \in \{-1,1\}$ be defined as in Theorem~\ref{th_NC_3eq}. There exists $T_*>0$ and $\varepsilon>0$ such that for any $T<T_*$, for every $s \in H^1((0,T),\R)$ with $s(0)=0$ and  $||\dot{s}||_{L^2} < \varepsilon$, the solution of system~(\ref{syst_aux}) satisfies
\begin{equation*}
\big( \tilde{\psi}^1(T), \tilde{\psi}^2(T) , \tilde{\psi}^3(T) \big) \neq e^{i \nu} \left( \Phi_1(T) e^{i \theta \mu} , \Phi_2(T) e^{i \theta \mu},
\left( \sqrt{1-\delta^2} + i \beta \delta \right) \Phi_3(T) e^{i \theta \mu} \right),
\end{equation*}
for all  $\delta >0 $, for all $\nu, \theta \in \R$.
\end{theo}

\noindent
The proof is very close to the one of Theorem~\ref{th_NC_aux_2eq_1}. 
\newline

\noindent
\textbf{Proof of Theorem~\ref{th_NC_aux_3eq} : }
Without loss of generality, we can assume $\BB>0$.
We consider the following quadratic form
\begin{align*}
\QQ_T(s) :&= \big( \lag \mu \varphi_3, \varphi_3 \rag - \lag \mu \varphi_2, \varphi_2 \rag \big)
\QQ_{T,1}(s)
+ \big( \lag \mu \varphi_1, \varphi_1 \rag - \lag \mu \varphi_3, \varphi_3 \rag \big)
\QQ_{T,2}(s)
\\
&+ \big( \lag \mu \varphi_2, \varphi_2 \rag - \lag \mu \varphi_1, \varphi_1 \rag \big)
\QQ_{T,3}(s),
\end{align*}
where $\QQ_{T,j}$ is defined as in (\ref{def_fq_js}).
This is rewritten as
\begin{equation*}
\QQ_T(s) = -\BB ||s||_{L^2}^2 + \int_0^T s(t) \int_0^t s(\tau) h(t-\tau) \md \tau \md t,
\end{equation*}
with $h \in C^0(\R,\R)$. 
Thus, there exists $T_*>0$, $C_* > 0$ such that for all $T<T_*$, 
\begin{equation*}
\QQ_T(s) \leq - C_* ||s||_{L^2}^2, \quad \forall s \in L^2((0,T),\R).
\end{equation*}

\noindent
Let $T<T_*$ and assume, by contradiction, that $\forall \varepsilon >0$, $\exists s_{\varepsilon} \in H^1((0,T),\R)$ with $s_\varepsilon(0)=0$ and $|| \dot{s}_{\varepsilon} ||_{L^2} < \varepsilon$ such that the associated solution of (\ref{syst_aux}) satisfies
\begin{equation*}
\big( \tilde{\psi}^1_{\varepsilon}(T) , \tilde{\psi}^2_{\varepsilon}(T) , \tilde{\psi}^3_{\varepsilon}(T) \big) = 
e^{i \nu_{\varepsilon}} \big( \Phi_1(T) e^{i \theta_{\varepsilon} \mu} , \Phi_2(T) e^{i \theta_{\varepsilon} \mu } ,
(\sqrt{1- \delta_{\varepsilon}^2} + i \beta \delta_{\varepsilon}) \Phi_3(T) e^{i \theta_{\varepsilon} \mu} \big),
\end{equation*}
with $\nu_{\varepsilon}, \theta_{\varepsilon} \in \R$ and $\delta_{\varepsilon} > 0$. Notice that,
\begin{equation*}
\delta_{\varepsilon} \underset{\varepsilon \to 0}{\rightarrow} 0, \quad
\theta_{\varepsilon} \underset{\varepsilon \to 0}{\rightarrow} 0, \quad
e^{i \nu_{\varepsilon}} \underset{\varepsilon \to 0}{\rightarrow} 1.
\end{equation*}

\begin{sloppypar} \noindent 
Straightforward computations and Lemma~\ref{lemme_approximation_aux} to estimate the terms $\Im \big( \lag \tilde{\psi}^1_{\varepsilon}(T) , \Phi_1(T) \rag \big) - \Im \big( \lag \tilde{\psi}^2_{\varepsilon}(T) , \Phi_2(T) \rag \big)$, 
$\Im \big( \lag \tilde{\psi}^1_{\varepsilon}(T) , \Phi_1(T) \rag \big)$ and $\Im \big( \lag \tilde{\psi}^3_{\varepsilon}(T) , \Phi_3(T) \rag \big)$ lead to
\begin{equation} \label{NC_estimation}
\theta_\varepsilon = \underset{\varepsilon \to 0}{o}(||s_\varepsilon||_{L^2}),
\quad
\sin(\nu_\varepsilon) = \underset{\varepsilon \to 0}{o}(||s_\varepsilon||_{L^2}),
\quad
\delta_\varepsilon = \underset{\varepsilon \to 0}{o}(||s_\varepsilon||_{L^2}).
\end{equation}
\end{sloppypar}


\noindent
For the sake of clarity, let us denote
\begin{align*}
\TT(s_{\varepsilon}):&= \big( \lag \mu \varphi_3, \varphi_3 \rag - \lag \mu \varphi_2, \varphi_2 \rag \big)
\Im(\lag \tilde{\psi}^1_{\varepsilon}(T),\Phi_1(T) \rag)
\\
&+ \big( \lag \mu \varphi_1, \varphi_1 \rag - \lag \mu \varphi_3, \varphi_3 \rag \big)
\Im(\lag \tilde{\psi}^2_{\varepsilon}(T),\Phi_2(T) \rag)
\\
&+ \big( \lag \mu \varphi_2, \varphi_2 \rag - \lag \mu \varphi_1, \varphi_1 \rag \big)
\Im(\lag \tilde{\psi}^3_{\varepsilon}(T),\Phi_3(T) \rag).
\end{align*}
Using estimates (\ref{NC_estimation}), straightforward computations lead to
\begin{equation*}
\TT(s_{\varepsilon}) = \beta \big( \lag \mu \varphi_2, \varphi_2 \rag - \lag \mu \varphi_1, \varphi_1 \rag \big) \cos(\nu_{\varepsilon}) \delta_{\varepsilon} + \underset{\varepsilon \to 0}{o}(|| s_{\varepsilon} ||_{L^2}^2).
\end{equation*}

\noindent
Finally, for $\varepsilon$ sufficiently small,
\begin{align*}
0 &< \beta \big( \lag \mu \varphi_2, \varphi_2 \rag - \lag \mu \varphi_1, \varphi_1 \rag \big) \cos(\nu_{\varepsilon}) \delta_{\varepsilon}
\\
&= \TT(s_{\varepsilon}) + \underset{\varepsilon \to 0}{o}(||s_{\varepsilon}||_{L^2}^2)
\\
&= \QQ_T(s_{\varepsilon}) + \underset{\varepsilon \to 0}{o}(||s_{\varepsilon}||_{L^2}^2)
\\
&\leq -C_* ||s_{\varepsilon}||_{L^2}^2 + \underset{\varepsilon \to 0}{o}(||s_{\varepsilon}||_{L^2}^2).
\end{align*}
This is impossible and ends the proof of Theorem~\ref{th_NC_aux_3eq}.

\hfill $\blacksquare$

\section{Conclusion, open problems and perspectives.}

\noindent
In this article, we have proved that the local exact controllability result of Beauchard and Laurent for a single bilinear Schr\"odinger equation cannot be adapted to a system of such equations with a single control. Thus, we developed a strategy based on Coron's return method to obtain controllability in arbitrary time up to a global phase or exactly up to a global delay for two equations. For three equations local controllability up to a global phase does not even hold in small time with small controls. Thus, in this setting and under generic assumptions no local controllability result can be proved in small time if $N \geq 3$. Finally, the main result of this article is the construction of a reference trajectory and application of the return method to prove local exact controllability up to a global phase and a global delay around $(\Phi_1,\Phi_2,\Phi_3)$.

\noindent
However our non controllability strategy is only valid for small time and we do not know if local exact controllability around the eigenstates $(\Phi_1,\Phi_2)$ hold in time large enough (for two equations or more). This would be the case if one manages to prove that the global delay $T^*$ can be designed to be the common period of the eigenstates $\Phi_k$ i.e. $T^* = \frac{2}{\pi}$. This is an open problem. 
Moreover, when Hypothesis~\ref{hypo_mu2} or \ref{hypo_mu3} are not satisfied, we do not know if the considered quadratic forms still have a sign. Thus, the question of non controllability when these hypotheses do not hold is an open problem. The question of non controllability with large controls has not been addressed here since our strategy relies on a second order approximation valid for small controls.

\noindent
The question of controllability of four equations or more is also open. In fact, each time we add an equation there is another diagonal coefficient $\lag \Psi^j, \Phi_j\rag$ which is lost. We proved that we can recover this lost direction using either a global phase or a global delay for $N=2$ and both a global phase and a global delay in the case $N=3$. It seems that there is no other degree of freedom to use to obtain controllability for $N \geq 4$. Moreover, there are other directions than the diagonal ones with the same gap frequencies (e.g. $\lambda_7 - \lambda_1 = \lambda_8 - \lambda_4$). Thus, for $N \geq 4$ one should consider a model with a potential that prevents such resonances.

\paragraph*{Acknowledgement} The author thanks K. Beauchard for having drawn his attention to this problem and for fruitful discussions.

\appendix
\section{Moment problems}
\setcounter{prop}{0}
\renewcommand{\theprop}{\Alph{section}\arabic{prop}}
\label{annexe_pb_moment}

We define the following space
\begin{equation*}
\ell^2_r (\N,\C) := \left\{ (d_k)_{k \in \N} \in \ell^2(\N,\C) \, ; \, d_0 \in \R \right\}.
\end{equation*}
In this article, we use several times the following moment problem result.
\begin{prop}
\label{prop_pb_moment}
Let T>0. Let $(\omega_n)_{n \in \N}$ be the increasing sequence defined by
\begin{equation*}
\left\{ \omega_n \, ; \, n \in \N \right\} = \left\{ \lambda_k - \lambda_j \, ; \, j \in \{1,2,3\}, \, k \geq j+1 \text{ and } k=j=3 \right\}.
\end{equation*}
There exists a continuous linear map
\begin{equation*}
\LL \, : \, \ell^2_r(\N,\C)  \rightarrow  L^2((0,T),\R),
\end{equation*}
such that for all $d :=(d_n)_{n\in \N} \in \ell^2_r(\N,\C)$,
\begin{equation*}
\int_0^T \LL(d)(t) e^{i \omega_n t} \md t = d_n, \quad \forall n \in \N.
\end{equation*}
\end{prop}

\noindent
\textbf{Proof of Proposition~\ref{prop_pb_moment} :} For $n \in \N^*$, let $\omega_{-n} := - \omega_n$. Using \cite[Theorems 9.1, 9.2]{KomornikLoreti05}, it comes that for any finite interval $I$, there exists $C_1, C_2 >0$, such that all finite sums
\begin{equation*}
f(t) := \sum_n c_n e^{i \omega_n t}, \quad c_n \in \C,
\end{equation*}
satisfy
\begin{equation*}
C_1 \sum_n |c_n|^2 \leq  \int_I |f(t)|^2 \md t \leq  C_2 \sum_n |c_n|^2.
\end{equation*}
This relies on Ingham inequality which holds true for any finite interval as $\lambda_k = k^2 \pi^2$.
Let $T>0$ and $H_0 := \text{Adh}_{L^2(0,T)} \big( \text{Span} \{ e^{i \omega_n \cdot } \, ; \, n \in \Z \} \big)$. Thus, $( e^{i \omega_n \cdot} )_{n \in \Z}$ is a Riesz basis of $H_0$ i.e.
\begin{equation*}
\begin{array}{cccc}
J_0: &  H_0
& \rightarrow &\ell^2(\Z,\C)
\\
& f & \mapsto & \left( \int_0^T f(t) e^{i \omega_n t} \md t \right)_{n \in \Z}
\end{array}
\end{equation*}
is an isomorphism (see e.g. \cite[Propositions 19, 20]{BeauchardLaurent}). Let $d \in \ell^2_r(\N,\C)$. We define $\tilde{d} := (\tilde{d}_n)_{n \in \Z} \in \ell^2(\Z,\C)$ by $\tilde{d}_n := d_n$, for $n \geq 0$ and $\tilde{d}_n := \overline{d_{-n}}$, for $n<0$. The map $\LL$ is defined by $\LL(d) := J_0^{-1} (\tilde{d})$. The construction of $\tilde{d}$ and the isomorphism property ensure that $\LL(d)$ is real valued.

\hfill $\blacksquare$

\bibliography{biblio} 
\bibliographystyle{plain} 

\end{document}